%% file: halfspace.tex
\documentclass[11pt]{article}

\usepackage{amsmath}
\usepackage{amsfonts}
\usepackage{amssymb}
\usepackage{amsthm}
\usepackage{graphics}
\usepackage{amscd}
\usepackage{graphicx,epsfig}
\usepackage{color}
\usepackage[all]{xy}
\usepackage{url}

\DeclareMathOperator{\Div}{div}

\renewcommand{\epsilon}{\varepsilon}
\newcommand{\boQ}{\mathcal{Q}}

\newcommand{\boD}{\mathcal{D}}

\newcommand{\boF}{\mathcal{F}}
\newcommand{\boC}{\mathcal{C}}

\newcommand{\boH}{\mathcal{H}}

\newcommand{\boO}{\mathcal{O}}
\newcommand{\boL}{\mathcal{L}}
\newcommand{\boS}{\mathcal{S}}

\newcommand{\boK}{\mathcal{K}}

\newcommand{\gotg}{\mathfrak{g}}

\newcommand{\gotf}{\mathfrak{f}}
\newcommand{\R}{\mathbb{R}}

\renewcommand{\H}{\mathbb{H}}

\newcommand{\N}{\mathbb{N}}
\newcommand{\la}{\langle}
\newcommand{\ra}{\rangle}
\newcommand{\eps}{\varepsilon}
\newcommand{\id}{\text{id}}
\newcommand{\der}[2]{\dfrac{\partial #1}{\partial #2}}
\newcommand{\dd}{\mathrm{d}}

\newcommand{\Ome}{\Omega}

\newtheorem{defn}{Definition}
\newtheorem{thm}{Theorem}

\newtheorem{prop}[thm]{Proposition}
\newtheorem{lem}[thm]{Lemma}

\renewcommand{\phi}{\varphi}
\newcommand{\dis}{\displaystyle}

\newtheorem*{thm*}{Theorem}

\newtheorem{claim}[thm]{Claim}
\newtheorem*{claim*}{Claim}

\theoremstyle{remark}

\newtheorem*{rem*}{Remark}
\newtheorem{example}{Example}

\newcounter{remark}

\newcounter{case}

\newcounter{construction}

\newcounter{fact}

\renewcommand{\d}{\mathbf{d}}
\newcommand{\Nil}{\mathrm{Nil}}
\newcommand{\Sol}{\mathrm{Sol}}

\title{A general halfspace theorem for constant mean curvature
surfaces}
\author{Laurent Mazet}

\begin{document}

\maketitle

\begin{abstract}
In this paper, we prove a general halfspace theorem for constant mean
curvature surfaces. Under certain hypotheses, we prove that, in an ambient
space $M^3$, any constant
mean curvature $H_0$ surface on one side of a constant mean curvature $H_0$
surface $\Sigma_0$ is an equidistant surface to $\Sigma_0$. The main
hypotheses of the theorem are that $\Sigma_0$ is parabolic and the mean
curvature of the equidistant surfaces to $\Sigma_0$ evolves in a
certain way.
\end{abstract}

\section{Introduction}

\input halfspace0.tex

\section{Preliminaries}\label{rappel}

\input halfspace1.tex

\section{Parabolic manifolds}
\label{secparabolic}

\input halfspace2.tex

\section{Regular $\eps$-neighborhood}
\label{neighborhood}

\input halfspace3.tex

\section{Construction of stable constant mean curvature surface}
\label{secstable}

\input halfspace4.tex

\section{The halfspace theorem}
\label{sechalfspace}

\input halfspace5.tex

\section{Halfspace theorems in certain ambient spaces}

\input halfspace6.tex

\noindent Universit\'e Paris-Est \\
Laboratoire d'Analyse et Math\'ematiques Appliqu\'ees, CNRS UMR8050\\
UFR des Sciences et Technologie\\ B\^atiment P3 4eme \'etage\\
61 avenue du G\'en\'eral de Gaulle\\
94010 Cr\'eteil cedex, France

\noindent \verb?laurent.mazet@math.cnrs.fr?

\end{document}

%% file: halfspace0.tex
One problem in the theory of constant mean curvature surfaces (cmc
surfaces) is to know
when two surfaces with the same constant mean curvature can coexist in the
same ambient space $M^3$. More precisely, if $\Sigma_1$ and $\Sigma_2$
are two properly immersed constant mean curvature $H_0$ surfaces in a
Riemannian $3$ manifold $M^3$ (these surfaces are called $H_0$ surfaces),
is the intersection $\Sigma_1\cap \Sigma_2$ empty?

If we consider two spheres in $\R^3$ with the same radius, we can put them
such a way that they do not meet. But inside a sphere of radius one, there
is no constant mean curvature one surface.

If we consider non intersecting properly immersed minimal surfaces in
$\R^3$, D.~Hoffman and
W.~Meeks \cite{HoMe} proved that these minimal surfaces are parallel
planes. For example, any minimal surface on one side of a plane is a plane.
This result is called a halfspace theorem.

This result can also be stated in an other way. Let us consider
a properly immersed minimal surface $\Sigma$ in $\R^3$ with compact
boundary and $P$ a plane. We assume that $\Sigma$ lies on one side of
$P$, then the distance between $\Sigma$ and $P$ satisfies
$d(P,\Sigma)=d(P,\partial\Sigma)$ \textit{i.e.} the distance is achieved
along the boundary. Such a result is called a maximum principle at
infinity. A very general maximum principle at infinity was proved by
W.~Meeks and H.~Rosenberg in \cite{MeRo}.

In a general setting, if $\Sigma_0$ is a properly embedded constant mean
curvature $H_0$ surface in $M^3$, a halfspace theorem with respect to
$\Sigma_0$ says that $H_0$ surfaces $\Sigma$ that lies on one side of
$\Sigma_0$ are ``classified''. Often, the classification implies that
$\Sigma$ has to be an equidistant surface to $\Sigma_0$. In this case, the
halfspace theorem can be interpreted as a maximum principle at infinity.

For example, A.~Ros and H.~Rosenberg \cite{RoRo} proved in $\R^3$
that no $H_0$ surface can lie in the mean convex side of a properly
embedded $H_0$ surface $\Sigma_0$ ($H_0>0$). We notice that this result
says that any $H_0$ surface in the mean convex side of $\Sigma_0$ is an
equidistant surface to $\Sigma_0$ but, since the equidistant surface to
$\Sigma_0$ do not have constant mean curvature $H_0$, no such surface can
exist.

Other halfspace theorems were proved by several authors. We have
halfspace theorems with respect to horospheres in $\H^3$ \cite{RodRos},
horocylinders in $\H^2\times\R$ \cite{HaRoSp}, vertical minimal planes in
$\Nil_3$ and $\Sol_3$ \cite{DaHa,DaMeRo} and entire minimal graph in
$\Nil_3$ \cite{DaMeRo}. We notice that, in \cite{DaMeRo}, B.~Daniel,
W.~Meeks
and H.~Rosenberg prove that the only minimal surfaces on one side of an
entire minimal graph in $\Nil_3$ are the vertical translate of the entire
graph. Since the distance between an entire graph and one of its translate
is not constant, the classification is of a different nature.

The aim of this paper is to give a general situation where a halfspace
theorem is true. More precisely, we prove that, under some
hypotheses, a $H_0$ surfaces that lies on one side of a given $H_0$ surface
is necessarily an equidistant surface.

Let $M^3$ be a complete Riemannian $3$ manifold which is geometrically 
bounded and $\Sigma_0$ a properly embedded constant mean curvature $H_0$
surface. Our main theorem says principally the following (see
Theorem~\ref{halfspace}, for a precise statement)

\begin{thm*}
Let $\Sigma_0\hookrightarrow M^3$ be as above. We assume that $\Sigma_0$
is parabolic.
\begin{enumerate}
\item Assume that the equidistant surfaces to $\Sigma_0$ has mean
curvature less than $H_0$ in the non mean convex side of $\Sigma_0$. Then
any $H_0$ surface that lies in the non mean convex side of $\Sigma_0$ and
is well oriented is an equidistant surface to $\Sigma_0$.
\item Assume that the equidistant surfaces to $\Sigma_0$ has mean
curvature larger than $H_0$ in the mean convex side of $\Sigma$. Then any
$H_0$ surface that lies in the mean convex side of $\Sigma_0$ is an
equidistant surface to $\Sigma_0$.
\end{enumerate}
\end{thm*}

In this result, the two important hypotheses are the parabolicity of
$\Sigma_0$ and the value of the mean curvature of the equidistant surfaces.
In fact, $\Sigma_0$ will be assumed to satisfy some other technical
hypotheses (see Theorem~\ref{halfspace}). 
The ``well oriented'' hypothesis means that, along the surface, the mean
curvature vector
points to $\Sigma_0$. When $\Sigma_0$ is a minimal surface ($H_0=0$), the
hypothesis on the mean curvature of the equidistant surfaces is that the
mean curvature vector does not point to $\Sigma_0$. In fact the hypothesis
about the mean curvature of the equidistant surface says that the mean
curvature evolves like the one of concentric spheres: inside the sphere of
radius $1$ the mean curvature is larger than $1$ outside it is less than
$1$.

If we consider $M^3=\R^3$ and $\Sigma_0$ is a plane. $\Sigma_0$ is
parabolic and the equidistant surface are also planes, thus the mean
curvature hypothesis is satisfied. The theorem then applies and we recover
the classical halfspace theorem.

Let us see why the hypotheses are important. We consider
$M^3=\H^2\times\R$ and the upper halfspace model for $\H^2$ \textit{i.e.}
$\H^2=\{(x,y)\in\R\times\R_+^*\}$ with the metric $\frac{1}{y^2} (\dd
x^2+\dd y^2)$. In
$\Ome=\{(x,y)\in\R\times\R^*,\ x>0\}$, we consider the function
$u(x,y)=\ln\frac{\sqrt{x^2+y^2}+y}{x}$. This function is a solution to
the minimal surface equation (its graph in $\H^2\times\R$ is a minimal
surface). As $x\rightarrow 0$, $u(x,y)\rightarrow +\infty$ and, as
$y\rightarrow 0$, $u(x,y)\rightarrow 0$. Let $\Sigma$ be the graph of $u$.
The minimal surface $\Sigma$ lies on one side of the minimal surface
$\Sigma_0=\H^2\times\{0\}$ and is asymptotic to it; so there is no
halfspace theorem for $\Sigma_0$. In fact the mean curvature of the
equidistant surfaces to $\Sigma_0$ is $0$. So the mean curvature 
hypothesis of the theorem is satisfied but $\Sigma_0$ is not parabolic.
The surface $\Sigma$ lies also on one side of the minimal surface
$\Sigma_1=\{x=0\}\times\R$. This
times, $\Sigma_1$ is parabolic (it is a flat $\R^2$) but the hypothesis for
the mean curvature of the equidistant surfaces is not satisfied. Thus, both
hypotheses are important in our statement.

Let us make a remark about the halfspace theorem of B.~Daniel,
W.~Meeks and H.~Rosenberg with respect to entire minimal graph in
$\Nil_3$ (Theorem~1.4 in \cite{DaMeRo}). Among all entire minimal graphs,
certain are not parabolic, so
their result is really of a different nature from the one we prove.

The paper is divided as follows. In the first section, we recall some
definition about constant mean curvature surfaces and we write the Stokes
formula in a general framework that we need.

In Section~\ref{secparabolic}, we explain what is a parabolic manifold and
we give a result that explain when the parabolicity is preserved by
quasi-isometry. In Section~\ref{neighborhood}, we explain what kind of
ambient space we consider in our halfspace theorem.

Section~\ref{secstable} is devoted to the proof of the first step of our
main theorem. It consist in proving that, if a $H_0$ surface lies on one
side of an other one, we can assume it is stable. In
section~\ref{sechalfspace}, we state our main theorem and finish its proof.

In the last section, we apply our main theorem to some ambient spaces that
have a Lie group structure. In this way, we recover known halfspace
theorems \cite{HoMe,RodRos,HaRoSp,DaHa,DaMeRo} and prove new results.

The author would like to thank H.~Rosenberg for many interesting and
helpful discussions

%% file: halfspace1.tex
In this section we recall some facts about cmc surfaces: what is the
stability and what can be said about self-intersection. We also explain
what is the Stokes formula in the setting of rectifiable boundary. Finally
we define the area estimate we will use in the following sections.

\subsection{Stability}

Let $S$ be a cmc surface in a Riemannian $3$-manifold $M$. On $S$, the
stability operator $L$ acts on smooth functions on $S$ by 
$$
Lu=-\Delta u-(2Ric(n,n)+|A|^2)u
$$
where $Ric(n,n)$ is the Ricci curvature of the ambient manifold, $n$ is the
unit normal vector to the surface and $|A|$ the norm of the second
fundamental form of $S$. $L$ is also called the Jacobi operator of $S$.

The cmc surface $S$ is said to be \emph{stable} if the stability operator
is non-negative on the set of smooth functions with compact support
\textit{i.e.}, for any smooth function $u$ with compact support,
$$
0\le\int_SuLu=\int_S\|\nabla u\|^2-(2Ric(n,n)+|A|^2)u^2.
$$
The stability operator appears as the second derivative of the area for
normal variations of the surface $S$ or as the first derivative of the mean
curvature (see \cite{BaDoCaEs}). 

\subsection{Stokes formula}\label{stokes}

Let $\Ome$ be a domain in $\R^n$ with a rectifiable boundary of finite
$\boH^{n-1}$ measure. This is the same as saying that the current
$[\Ome]$ associated to $\Ome$ has a rectifiable boundary
$\partial[\Ome]$. By Theorem 4.1.28 and Theorem 4.5.6 in \cite{Fed} (see
also 4.5.12 in \cite{Fed} and 12.2 in \cite{Mor}), for any smooth vector
field $X$ with compact support in $\R^n$, the Stokes formula can be
written:
\begin{equation}\label{forstokes}
\int_{\Ome}\Div X(x)\dd\boL^nx=\int_{\partial\Ome}X(x)\cdot
n(\Ome,x)\dd\boH^{n-1}x
\end{equation}
where $n(\Ome,x)$ is a unit vector called the exterior normal of $\Ome$ at
$x$ (see 4.5.5 in \cite{Fed} for a definition in this situation). This
exterior normal is defined $\boH^{n-1}$ almost everywhere along
$\partial\Ome$. We notice that the definition of $n(\Ome,x)$ is local and
coincides with the classical unit outgoing normal vector for smooth
boundaries.

Moreover, the $\boH^{n-1}$ measure of $\partial\Ome$ is equal to the mass
of the $n-1$ current $\partial[\Ome]$.

\subsection{Self intersection}\label{selfinter}

Now let us consider $D_2$, an open disk in $\R^2$, and $D_1$, the open disk
with the same center and half radius. On $D_2\times\R$, we consider a
Riemannian metric $g$. Let $f_1,\cdots,f_n$ be smooth functions on
$D_2$ such that their graphs have constant mean curvature $H_0$ with
respect to the metric $g$ and the mean curvature vector points downward.
Let $p$ in $D_2$ such that $f_i(p)=f_j(p)$ and $\nabla(f_i-f_j)(p)=0$, $p$
is a singular intersection point. The structure of the set $\{f_i=f_j\}$
near $p$ is then described by Theorem 5.3 in \cite{CoMi2}: it is the union
of $2d$ embedded arcs meeting at $p$. Moreover such points are
isolated.

Let $f_0$ be a smooth function on $D_2$ such that $\nabla(f_0-f_i)$ never
vanishes wherever $f_0=f_i$. We notice that if $p$ satisfies
$f_i(p)=f_j(p)$ and $\nabla(f_i-f_j)(p)\neq0$, the level set $\{f_i=f_j\}$
is locally an embedded arc. This implies that $I_{i,j}=\{f_i=f_j\}$ is
locally either a smooth arc or the union of embedded arcs meeting at a
point. Thus $I_{i,j}\cap D_1$ is compact.

We define the function $f$ by $f(p)=\min_if_i(p)$. Let $\Ome_i$ be the
open subset $\Ome_i=\{p\in D_1|\, f(p)=f_i(p) \textrm{ and }\forall j\neq
i\ f(p)<f_j(p)\}$. The question is: what is the Stokes formula for such a
domain $\Ome_i$?

First we see that $\partial\Ome_i$ is included in the sets $I_{i,j}$ and
the boundary of $D_1$. This implies that $\partial\Ome_i$ is a
$1$-rectifiable subset of finite $\boH^1$ measure. Thus the above
formula~\eqref{forstokes} can be applied. But we need to understand what is
$n(\Ome_i,x)$ and
where it is defined.

Let $p$ be a point in $\partial\Ome_i\setminus \partial D_1$. First, since
 the singular intersection points form a discrete set, this set is finite
in $D_1$ and has vanishing $\boH^1$ measure. So we assume that $p$ is not
such a point. We denote by $A(p)$ the set of indices $j$ such that $p\in
\partial\Ome_j$; we have $i\in A(p)$ and for any $j\in A(p)$ $f(p)=f_j(p)$.
There is two situations.

First, the vectors $\nabla(f_j-f_l)(p)$ for $j\neq l$ and
$j,l\in A(p)$ are not all linearly dependant (this implies that the
intersection of the tangent planes to the graphs of the $f_j$, $j\in A(p)$,
is a point). In this case, near $p$ the domain $\Ome_i$ is included in an
angular sector of angle strictly less than $\pi$. This implies that the
exterior normal $n(\Ome_i,p)$ does not exist. This is the same for
$n(\Ome_j,p)$, $j\in A(p)$.

Let us assume now that the vectors $\nabla(f_j-f_l)(p)$ for $j\neq l$ and
$j,l\in A(p)$ are all linearly dependant: the intersection of the tangent
planes to the graphs is now a line (this is the case when $A(p)$ has only
two elements). In this case, all the curves $I_{j,l}$, $j,l\in A(p)$, are
tangent at $p$. Let $L_j$ be the differential of $f_j$ at $p$. For any
$x\in\R^2$, we define $L(x)=\min_{j\in A(p)}L_j(x)$. Since the $L_j$ are
linear and different, there exists, in fact, a unique subset
$\{j_1,j_2\}\subset A(p)$ such that $L(x)=\min(L_{j_1}(x),L_{j_2}(x))$.
$\{L=L_{j_1}\}$ and $\{L=L_{j_2}\}$ are half-planes and we denote by $n$
the unit vector normal to $\{L_{j_1}=L_{j_2}\}$ and pointing in
$\{L=L_{j_2}\}$. Thus, for any $\lambda<1$, the affine angular
sector $\{p+x, \lambda\|x\|<n\cdot x\}$ is included in $\Ome_{j_2}$ near
$p$ and the affine angular sector $\{p+x, \lambda\|x\|<-n\cdot x\}$ is
included in $\Ome_{j_1}$ near $p$. This implies that $n(\Ome_{j_1},p)=n$,
$n(\Ome_{j_2},p)=-n$ and $n(\Ome_{j},p)$ is not defined for any $j\in
A(p)\setminus \{j_1,j_2\}$. In this case we say that $p$ is in the set
$\Gamma_{j_1,j_2}=\Gamma_{j_2,j_1}$. 

Finally, with this definition, if $X$ is a smooth vector field with
compact support in $D_1$, we get the following Stokes formula:
$$
\int_{\Ome_i}\Div X(x)\dd\boL^2x=\sum_{j\neq i}
\int_{\Gamma_{i,j}}X(x)\cdot n(\Ome_i,x)\dd\boH^1x.
$$

\subsection{Area bounds}\label{areabound}

In this subsection, we define a notion of area bound. Let $M$ be a
Riemannian $3$-manifold. Let $p$ be a point in $M$ and $P$ be a plane in
$T_pM$. Let $(e_1,e_2,e_3)$ be an orthonormal basis of $T_pM$ such that
$P$ is the plane generated by $e_1$ and $e_2$. We denote by $E_{p,P}(t)$
be the image by the exponential map at $p$ of the ellipsoid
$\{(x,y,z)\in T_pM|\, x^2+y^2+4z^2\le t^2\}$ where $(x,y,z)$ are the
coordinates in $T_pM$ with respect to $(e_1,e_2,e_3)$.

Let $V(O_{p,P}(t))$ be the volume of $O_{p,P}(t)$ and
$A(\partial O_{p,P}(t))$ be the area of its boundary. We notice that, for
$t$ small, $V(O_{p,P}(t))\sim (2/3)\pi t^3$ and $A(\partial O_{p,P}(t))\sim
4\pi\alpha t^2$ with $\alpha <1$.

\begin{defn}
Let $(S_i)_i$ be a family of immersed surfaces in a Riemannian
$3$-manifold $M$. We say that the family satisfies a uniform area
estimate of at most one leaf if for any $p\in M$ and $P$ a plane in
$T_pM$
there exists $t_0>0$ and $\beta<1$ such that,  for any $t<t_0$ and $i$,
$$A(S_i\cap O_{p,P}(t))\le 2\beta\pi t^2.$$
\end{defn}

If $S$ is an immersed surface and $p$ is a point in $S$, we have
$$
\liminf_{t\rightarrow 0}\frac{A(S\cap O_{p,T_pS}(t))}{t^2}\ge \pi.
$$
Thus the area estimate of at most one leaf prevents $S$ to
pass at $p$ more than one time with the same tangent plane.

%% file: halfspace2.tex
In this section, we recall some definitions about the conformal type of
Riemannian manifolds and we explain when the conformal type is
preserved by quasi-isometry. We refer to \cite{Gri} for a general
presentation of conformal types.

Let $(M,g)$ be a Riemannian manifold. A continuous function $u$ on a domain
$\Ome\in M$ is \emph{superharmonic} if, for any precompact domain
$U\subset\subset\Ome$ and any harmonic function $v\in C^2(U)\cap
C^0(\overline{U})$, $v\le u$ on $\partial U$ implies $v\le u$ on $U$. If
$u_1,\cdots,u_n$ are superharmonic functions, we remark that $u=\inf_i u_i$
is also a superharmonic function. 

\begin{defn}
Let $(M,g)$ be a  Riemannian manifold.
\begin{enumerate}
\item If $\partial M=\emptyset$, $M$ is \emph{parabolic} if any bounded
superharmonic function on $M$ is constant.
\item If $\partial M\neq \emptyset$, $M$ is \emph{parabolic at infinity} if
any bounded non-positive superharmonic function on $M$ with $u=0$ on
$\partial M$ is constant.
\end{enumerate}
\end{defn}

When $\partial M\neq\emptyset$, $M$ is often said to be ``parabolic" instead
of ``parabolic at infinity", but we prefer to use different terminologies.
In fact, there is a lot of equivalent characterizations of parabolicity
(see \cite{Gri}) and we will use certain of them below. As an example, a
Riemannian manifold $M$  without boundary is parabolic if and only if there
exists a sequence $(\phi_n)_n$ of smooth functions with compact support in
$M$ such that $0\le \phi_n\le 1$, $(\phi_n^{-1}(1))_n$ is an
increasing exhaustion by compact subsets of $M$ and
$$
\lim_{n\rightarrow +\infty}\int_M\|\nabla\phi_n\|^2=0
$$
We remark that
a subdomain of a parabolic manifold, viewed as a manifold with
boundary, is parabolic at infinity.

Let $(M,g)$ and $(N,h)$ be two $n$-dimensional Riemannian manifold and
let $F$ be a map from $M$ to $N$. If $k\ge 1$, we say that $F$ is $k$
quasi-isometric or a local $k$ quasi-isometry if, for any $p\in M$ and
$v\in T_p M$, we have $\frac{1}{k}\|v\|_g\le \|T_p F(v)\|_h\le k\| v\|_g$.
If
$M$ and $N$ has no boundary and $F:M\rightarrow N$ is a $k$
quasi-isometric diffeomorphism, $M$ is
parabolic if an only if $N$ is parabolic. For parabolicity at infinity, we
do not have such a result. In fact we have the following proposition:

\begin{prop}\label{parabolic}
Let $(M,g)$ and $(N,h)$ be two $n$-dimensional Riemannian manifold such
that $\partial M\neq\emptyset$ and $N$ has no boundary. We assume that
$(N,h)$ is parabolic and that there exists $F:M\rightarrow N$ an injective
local $k$ quasi-isometry. Then $M$ is parabolic at infinity.
\end{prop}

\begin{proof}
Let us assume that $M$ is not parabolic at infinity, then it exists a
harmonic function $u_M$ such that $0<u_M\le 1$, $u_M=1$ on $\partial M$ and
$\inf_M u_M=0$ (see \cite{Gri}, $u_M(x)$ is the probability that a
Brownian motion from $x$ hits the boundary of $M$). Let $\eta\in(0,1)$ be
a regular value of $u_M$ and $\phi\in C^\infty(\R,[0,1])$ be a function
such that $\phi=0$ on $[(1+\eta)/2,+\infty)$ and $\phi=1$ on
$(-\infty,\eta]$.

Since $F$ is a quasi-isometry, $DF$ is invertible and
$\tilde{g}=F_*(g)$ is well defined. $\tilde{g}$ is a section over $F(M)$
of the symmetric 2-tensor bundle. Moreover $\tilde{g}$ is positive definite
and we have $\frac{1}{k^2}h\le \tilde{g}\le k^2 h$.

We denote by $v$ the function $u_M\circ F^{-1}$ on $F(M)$ and we consider
$\mu=\phi\circ v$. The function $\mu$ is $C^\infty$ on $F(M)$ and
vanishes on $v^{-1}([(1+\eta)/2,1])$. This domain contains a neighborhood
of $F(\partial M)=\partial F(M)$. So we can extend the definition of the
function $\mu$ by $0$ to the complement of $F(M)$. The function $\mu$ is then 
$C^\infty$ on $N$ with $\mu=1$ on $v^{-1}([0,\eta])$.

On $N$, we define $\tilde{h}=(1-\mu)h+ \mu \tilde{g}$ ($\mu
\tilde{g}$ is well defined on $N$ since $\mu$ vanishes outside $F(M)$).
$\tilde{h}$ is a global section of the symmetric 2-tensor bundle and we
have:
$$
\frac{1}{k^2}h\le ((1-\mu)+\frac{1}{k^2}\mu)h\le \tilde{h}\le
((1-\mu)+k^2\mu)h\le k^2 h
$$ 
So $\tilde{h}$ defines a Riemannian metric on $N$
and $\id_N:(N,h)\rightarrow (N,\tilde{h})$ is a local quasi-isometry.

Since $(N,h)$ is parabolic, so is $(N,\tilde{h})$. Let $\tilde{v}$ be the
function defined by $\eta$ outside $F(M)$ and by $\min(\eta,v)$ on $F(M)$;
$\tilde{v}$ is continuous on $M$. On $v^{-1}([0,\eta])$,
$\tilde{h}=\tilde{g}$ so $v^{-1}([0,\eta])$ with the metric $\tilde{h}$ is
isometric by $F$ with $u_M([0,\eta])\subset M$. Thus, on
$v^{-1}([0,\eta])$,
$\Delta_{\tilde{h}}\tilde{v}=\Delta_{\tilde{g}}v=(\Delta_g u_M)\circ
F^{-1}=0$ since $u_M$ is harmonic. On the complement to
$v^{-1}([0,\eta])$, $\tilde{v}$ is constant so
$\Delta_{\tilde{h}}\tilde{v}=0$. So $\tilde{v}$ is a positive
superharmonic function on $(N,\tilde{h})$ (it is locally the infimum of
two harmonic function) and it is bounded from above by
$\eta$. $\tilde{v}$ is then constant and equal to $\eta$. This implies
that $u_M=\eta$ on $u_M^{-1}([0,\eta])$ which contradicts $\inf_M u_M=0$.
This ends the proof of the proposition.
\end{proof}

%% file: halfspace3.tex
In this section, we explain what kind of ambient space we will consider
in our main theorem. Let $\Sigma$ be a properly embedded
constant mean curvature $H_0$ surface in an ambient $3$-manifold
$M$. The
$\eps$-tubular neighborhood of $\Sigma$ is the set of points in $M$ at
distance less than $\eps$ from $\Sigma$. We can define the map
$F:\Sigma\times[-\eps,\eps]\rightarrow M,\ (x,t)\mapsto \exp_x(tn(x))$
where $n(p)$ is the unit normal vector such that the mean curvature vector
of $\Sigma$ at $p$ is $-H_0n(p)$. The image of $F$ is the $\eps$-tubular
neighborhood of $\Sigma$. When $F$ is a diffeomorphism, it gives a global
parametrization of the neighborhood. Besides if $H_0>0$, the image of
$F(\Sigma\times[-\eps,0])$ is the mean convex side of the tubular
neighborhood and $F(\Sigma\times[0,\eps])$ is the non-mean convex side.
When $H_0=0$, no such distinction can be done.

We want to take this situation as a model for our ambient spaces.
\begin{defn}
Let $(\Sigma,\dd\sigma_0^2)$ be a $2$-dimensional complete Riemannian
manifold. An \emph{outside $\eps$-half neighborhood} of $\Sigma$ is the
$3$-manifold with boundary $M_+(\eps)=\Sigma\times[0,\eps]$ with a
Riemannian metric $\dd s^2=\dd \sigma_t^2+\dd t^2$ where $t\mapsto
\dd\sigma_t^2$ is a smooth family of Riemannian metric on $\Sigma$ such
that $\dd s^2$ is complete. 

An \emph{inside $\eps$-half neighborhood} of $\Sigma$ is the
$3$-manifold with boundary $M_-(\eps)=\Sigma\times[-\eps,0]$ with a
Riemannian metric $\dd s^2=\dd \sigma_t^2+\dd t^2$ where $t\mapsto
\dd\sigma_t^2$ is a smooth family of Riemannian metric on $\Sigma$ such
that $\dd s^2$ is complete. 
\end{defn}

It seems that we define two times the same object but we prefer to use two
different terms for the model of the mean convex side (the inside
$\eps$-half neighborhood) and the non-mean convex side (the outside
$\eps$-half neighborhood).

Let $M_\pm(\eps)$ be a $\eps$-half neighborhood of $\Sigma$. If $\eps'\le
\eps$, the submanifold $\Sigma\times[0,\eps']\subset M_+(\eps)$ is denoted
by $M_+(\eps')$ and is an outside $\eps'$-half neighborhood.
$M_-(\eps')=\Sigma\times[-\eps',0]\subset M_-(\eps)$ is an inside
$\eps'$-half neighborhood.

We denote by $\Sigma_t$ the submanifold $\Sigma\times\{t\}$, $\Sigma_0$
with its induced metric is then isometric to $(\Sigma,\dd \sigma_0^2)$.
$\Sigma$ is then isometrically embedded in $M_\pm(\eps)$. $\Sigma$ and
 $\Sigma_0$ will be often viewed as the same object. We denote
$M_\pm^*(\eps)=M_\pm(\eps)\setminus \Sigma_0$. We
also define the distance function $\d$ as $\d(x,t)=|t|$, $\d$ is then the
distance from $\Sigma_0$. $\Sigma_t$ is the equidistant surface from
$\Sigma_0$ at distance $|t|$. On $M_\pm(\eps)$, we define the projection
map
$\pi:M_\pm(\eps)\rightarrow \Sigma_0$ by $\pi(x,t)=(x,0)$. We denote by
$\pi_t$ the restriction of $\pi$ to $\Sigma_t$.

Let $\xi$ denote the unit vectorfield $\der{}{t}$. In
the following, we always consider $-\xi$ as the unit normal vector to the
surface $\Sigma_t$. This is with respect to this normal vector that we will
compute the mean curvature. For $(x,t)\in M_\pm(\eps)$, this implies that
$\Div 
\xi(x,t)=2H(x,t)$ where $H(x,t)$ is the mean curvature of $\Sigma_t$ at
$(x,t)$.

Let $H_0\ge 0$ be a constant we say that $M_+(\eps)$ satisfies the
\emph{$H\le H_0$ hypothesis} if:
\begin{center}
for any $t\in[0,\eps]$ the mean curvature of $\Sigma_t$ is less than
$H_0$ at any point. 
\end{center}
We say that $M_-(\eps)$ satisfies the \emph{$H\ge H_0$ hypothesis} if:
\begin{center}
for any $t\in[-\eps,0]$ the mean curvature of $\Sigma_t$ is larger than
$H_0$ at any point. 
\end{center}

\begin{defn}
Let $(\Sigma,\dd \sigma_0^2)$ be a complete $2$-dimensional Riemannian
manifold. Let $M_\pm(\eps)$ be an outside or inside
$\eps$-half neighborhood. We say that $M_\pm(\eps)$ is \emph{regular} if
\begin{enumerate}
\item there is $k>0$ such that $\pi_t$ is a $k$ quasi-isometric map for any
t with $|t|\le \eps$.
\item there is $C$ such that the norm of the second fundamental form of
$\Sigma_t$ is bounded by $C$ for any t with $|t|\le \eps$.
\item $M_\pm(\eps)$ is geometrically bounded.
\end{enumerate}
\end{defn}

Let $S$ be a properly immersed cmc $H_0$ surface ($H_0\ge 0$) in $M_{\pm}(\eps)$ with
$S\subset M_\pm^*(\eps)$ and 
possibly nonempty boundary in $\Sigma_{\pm\eps}$. Along $S$, we always consider the unit normal vector $n$ such that the mean curvature vector is $H_0 n$. We denote by $D$ the
connected component of $M_\pm(\eps)\setminus S$ that contains $\Sigma_0$.
Consider a point $p$ in $S\cap \partial D$. Let $\Delta\subset S$ an embedded neighborhood  of $p$. Let us
consider the map $F:\Delta\times(-\eta,\eta)\rightarrow M_\pm(\eps),
(q,t)\mapsto \exp_q(tn(q))$. $F$ is an embedding if $\eta$ is small
enough and its image is a neighborhood of $p$ in $M_\pm(\eps)$. We say that
the mean curvature vector of $S$ at $p$ points into $D$ (resp. not into
$D$) if, for any sequence $(p_n)_n$ in $D$ with $p_n\rightarrow p$, $p_n\in
F(\Delta\times [0,\eta))$ (resp $p_n\in F(\Delta\times (-\eta,0])$) for large $n$.

We say that $S$ is \emph{well oriented} if, for any point in $S\cap
\partial D$, the mean curvature vector of $S$ points into $D$ (resp. not
into $D$) when $S\looparrowright M_+(\eps)$ (resp. $S\looparrowright
M_-(\eps)$). We notice that, when $S$ is minimal ($H_0=0$), $S$ can be assumed to be orientable by considering a covering space. Moreover it can always be considered as well oriented.

%% file: halfspace4.tex
Let $M_\pm(\eps)$ be an $\eps$-half neighborhood of a surface
$(\Sigma,\dd\sigma_0^2)$ such that $\Sigma_0$ has constant mean curvature
$H_0$. The main result of our paper says under which hypotheses we have a
halfspace theorem for $\Sigma_0$ : any properly immersed constant
mean curvature $H_0$ surface in $M_\pm(\eps)$ is an equidistant surface to
$\Sigma_0$. In this section, we explain that, if such a constant mean
curvature $H_0$ surface exists, we can assume that it is stable.

\begin{thm}\label{constab}
Let $(\Sigma_0,\dd \sigma_0^2)$ be a complete orientable Riemannian
surface, $\eps$
be positive and $H_0$ be  non-negative. Let $M_\pm(\eps)$ be an inside or
outside $\eps$-half neighborhood of $\Sigma$. We consider a
properly immersed constant mean curvature $H_0$ surface $S$ in
$M_\pm(\eps)$
with possibly nonempty boundary in $\Sigma_\eps$ and $S\subset
M_\pm^*(\eps)$. We assume that the
lower bound of the distance function $\d$ on $S$ is $0$.
\begin{enumerate}
\item If $S\looparrowright M_+(\eps)$ is well oriented and $M_+(\eps)$
satisfies the $H\le H_0$ hypothesis. There exist $\eps'>0$ and a
properly immersed constant mean curvature $H_0$ surface $S'$ in
$M_+(\eps')$ with non empty boundary in $\Sigma_{\eps'}$ such that
$S'\subset M_+^*(\eps')$, $S'$ is stable, well
oriented and the distance function $\d$ on $S'$ is not constant.
\item If $S\looparrowright M_-(\eps)$ and $M_-(\eps)$ satisfies the $H\ge
H_0$ hypothesis. There exist $\eps'>0$ and a properly immersed constant
mean curvature $H_0$ surface $S'$ in $M_-(\eps')$
with non empty boundary in $\Sigma_{\eps'}$ such that $S'\subset M_-^*(\
\eps')$, $S'$ is stable and
the distance function $\d$ on $S'$ is not constant.
\end{enumerate}
\end{thm}

The remaining part of the section is devoted to the proof of this result.
But let us begin by some remarks on the proof and the result.

The first remark is that, for $H_0=0$, both cases are in fact the same 
since the good orientation hypothesis has no meaning and the outside and
inside half neighborhoods play the same role for minimal surfaces.

Besides the proof of both cases are very similar, so, the detailed proof is
written for the first case with $H_0>0$. Then we explain the differences
for the other ones.

One more remark about this result is that, if $S$ is stable or if
a surface $S\cap M_\pm(\eps')$, for $\eps'<\eps$, is stable then this
surface gives the surface $S'$ we look for. If no such surface is stable,
the surface $S'$ produced by the proof is, in fact, embedded and well
oriented in both cases.

A large part of the proof is inspired by the work of A.~Ros and
H.~Rosenberg in \cite{RoRo} and L.~Hauswirth, P.~Roitman and H.~Rosenberg
in \cite{HaRoRo}.

\subsection{$S\looparrowright M_+(\eps)$ and $H_0>0$}

Let us consider $S$ in $M_+^*(\eps)$ with $H_0>0$. First we need
to introduce objects that will be used in the proof.

\subsubsection{Definition of the barriers}

Let $x_0$ be a point in $\Sigma$ and $\eta_0>0$ such that the exponential
map $\exp_{x_0}$ for the metric $\dd\sigma_0^2$ is a diffeomorphism from
the disk of radius $\eta_0$ in $T_{x_0}\Sigma$ into a neighborhood
$D_{\eta_0}$ of $x_0$. Since $S$ is properly immersed, there is
$\eps_0$ such that $D_{\eta_0}\times[0,\eps_0]\cap S=\emptyset$. In
$D_{\eta_0}\times[0,\eps_0]$ we consider the chart $\exp_{x_0}\times\id$
defined on $\Delta_{\eta_0}\times[0,\eps_0]$ where $\Delta_{\eta_0}$
is the Euclidean disk in $T_{x_0}\Sigma$ of radius $\eta_0$. Let $\eta$ be
small and consider in $\Delta_{\eta_0}\times[0,\eps_0]$ the surfaces of
revolution $C_{\eta,t}$ parametrized by
$$
X_{\eta,t}(u,v)=\left((t-\frac{\eta}{6}\cos v)\cos u,(t-\frac{\eta}{6}\cos
v)\sin u,\frac{\eta}{6}(1+\sin v)\right)
$$
where $(u,v)\in[0,2\pi]\times[-\pi/2,\pi/2]$, $t\le \eta_0$ and $\
\eta\le \min(\eta_0, \eps_0)$.

Let $\eta$ be sufficiently small so that the surfaces $C_{\eta,t}$ are
well defined for $t\in[\eta/2,\eta]$. We denote by $K$ the compact domain
of $\Delta_{\eta_0}\times[0,\eta/3]$ bounded by $C_{\eta,\eta/2}$ and
containing the origin. For $t\in[\eta/2,\eta]$, we denote by $Q_t$ the
domain of $\Delta_{\eta_0} \times[0,\eta/3]$ bounded by $C_{\eta,\eta/2}$
and $C_{\eta,t}$ ($Q_t\subset Q_\eta$). $Q_\eta$ is foliated by the
surfaces $C_{\eta,t}$ for $t\in[\eta/2,\eta]$. On these
surfaces, the mean curvature vector does not point to $K$ and its
norm is larger than $1/\eta$. We denote by $\boK_{bar}$, $\boQ_{bar,t}$ and
$\boC_{\eta,t}$ the images of $K$, $Q_t$ and $C_{\eta,t}$ in
$D_{\eta_0}\times[0,\eps_0]\subset M_+(\eps_0)$ (see Figure~\ref{kbar}),
we also denote $\overset{\circ}{\boQ_{bar,t}}=\boQ_{bar,t}\setminus
\boC_{\eta,t}$. In
fact, $\boK_{bar}\cup \boQ_{bar,\eta}\subset M_+(\eta/3)$ and these subsets
do not meet $S$. Let $\xi_{bar}$ be the unit vector normal to
$\boC_{\eta,t}$ which does not point to $\boK_{bar}$. $\xi_{bar}$ is a unit
vector field on
$\boQ_{bar,\eta}$ and, because of the value of the Euclidean mean
curvature, we can choose $\eta$ sufficiently small such that $\Div
\xi_{bar}$ is as small as we want. So we choose $\eta$ such that $\Div
\xi_{bar}\le 2H_0$ and $\eta/3$ is a regular value of the function $\d$ on
$S$.

We write $\eps_1=\eta/3$. From now on, we work in $M_+(\eps_1)$ and we
consider the restriction of $S$ to $M_+(\eps_1)$ which we still call $S$.

\begin{figure}[h]
\begin{center}
\resizebox{1\linewidth}{!}{\input{fighalfspace1.pstex_t}}
\caption{}
\label{kbar}
\end{center}
\end{figure}


As explained above, if $S$ is stable  Theorem~\ref{constab} is already
proved so we can assume $S$ is not stable. Hence there exists an
exhaustion $(K_n)_n$ of $\Sigma$ by compact subsets such that, for any $n$,
$S_n=S\cap K_n\times[0,\eps_1]$ is not stable.

We denote by $D$ the connected component of $M_+(\eps_1)\backslash S$
that contains $\Sigma_0$ and $D_n=D\cap
(K_n\times[0,\eps_1])$.
We notice that $\boK_{bar}\cup \boQ_{bar,\eta}\subset D_n$ for large $n$.
Since $S$ is well oriented, the mean curvature vector of $S$ points into
$D_n$ along $S\cap \partial D_n$.
Let $\phi$ be the first eigenfunction of the Jacobi operator of $S_n$ :
$\phi$ vanishes on $\partial S_n$, is positive at the interior of $S_n$ 
and satisfies $-L\phi+\lambda_1\phi=0$ where $L$ is the stability operator
and $\lambda_1$ is a negative constant. Perturbing the $K_n$, we can
assume that $0$ is not an eigenvalue of $-L$, hence  there is a
smooth function $v$ on $S_n$, vanishing on the boundary such that $-Lv=1$
in $S_n$. By the boundary maximum principle, the outing derivative
$\der{\phi}{\nu}$ is negative along $\partial S_n$. Thus for $a$ small
enough, the function $u=\phi+av$ is positive in the interior of $S_n$. 

Let $N(x)$ be the unit normal to $S_n$ such the mean curvature vector is
$H_0N(x)$. For $t_0>0$, we define, on $S_n\times[0,t_0]$, the map
$F(x,t)=\exp_x(tu(x)N(x))$, we assume $t_0$ small such that $F$ is
an immersion. We then denote $\widetilde{\boQ}_{uns}^n=S_n\times[0,t_0]$
with the induced metric $F^*\dd s^2$. In $\widetilde{\boQ}_{uns}^n$ we
consider the surfaces $\boS_t^n=S_n\times\{t\}$ which foliates
$\widetilde{\boQ}_{uns}^n$. Let $\xi_{uns}$ be the unit vector
field defined on $\widetilde{\boQ}_{uns}^n$ normal to $\boS_t^n$. We have
$\Div \xi_{uns}=-2H_t$ where $H_t$ is the mean
curvature of $\boS_t^n$. Moreover, we have:
$$
\frac{d}{dt}_{|t=0}2H_t=-L'u=-\lambda_1\phi+a>0
$$
so, choosing $t_0$ small enough, we get $H_t>H_0$. We define
$\boQ_{uns}^{n,0}=F(\widetilde{\boQ}_{uns}^n)\cap D_n$ and
$\boD_n=D_n\setminus \boQ_{uns}^{n,0}$ (see Figure~\ref{quns}). In fact,
$t_0$ is also chosen such that $\boD_{n-1}\cap (K_{n-2}\times[0,\eps_1])
\subset \boD_n\cap (K_{n-2}\times[0,\eps_1])$. This implies that the
sequence $(\boD_n)_n$ is increasing with respect to compact subsets. We
can also assume that $\cup_n\boD_n=D$.

The surface $\boS_{t_0}^n$ is immersed by $F$ in $M_+(\eps_1)$ and the
normal vector $F_*(\xi_{uns})$ points to $\boD_n$ along $\boS_n$ where
$\boS_n=F(\boS_{t_0}^n)\cap \partial \boD_n$. Let $x$ be a
point in $\boS_{t_0}^n$ and consider
$D_x\subset \boS_{t_0}^n$ a small open geodesic disk which is embedded in
$M_+(\eps_1)$ by $F$. Let $\psi$ be a smooth function on $\boS_{t_0}^n$
vanishing outside $D_x$ and positive in $D_x$. We then define on
$\boS_{t_0}^n\times[0,2t_x]$ 
$$
G(p,t)=\exp_{F(p)}(t\psi(p)F_*(\xi_{uns}(p)))
$$
If we choose $t_x$ small enough, we can assume that $G$ is an embedding on
$D_x\times[0,2t_x]$. In $G(D_x\times[0,2t_x])$, we define $\xi_x$ the unit
vector
field normal to the embedded surfaces $\boS_t^x=G(D_x\times\{t\})$ with
$\xi_x=F_*(\xi_{uns})$ along $\boS_0^x$. Since the mean curvature of
$\boS_0^x$ is larger than $2H_0$, if $t_x$ is small enough, we can assume
that $\Div \xi_x<-2H_0$.

For $\delta\in[1,2]$, we denote $\boQ_{x,\delta}=G(D_x\times[0,\delta
t_x])$ (see Figure~\ref{quns}) and
$\overset{\circ}{\boQ_{x,\delta}}=G(D_x\times[0,\delta
t_x))$. Now we define
$$
\boQ_{uns}^{n,1}=\boD_n\cap\bigcup_{x\in
\boS_{t_0}^n}\boQ_{x,1}.
$$

Since $F_*(\xi_{uns})$ points to $\boD_n$ along $\boS_n$, any point in
$\boS_n$ is at a positive distance from
$\boD_n\setminus \boQ_{uns}^{n,1}$.

\begin{figure}[h]
\begin{center}
\resizebox{1\linewidth}{!}{\input{fighalfspace2.pstex_t}}
\caption{}
\label{quns}
\end{center}
\end{figure}


Let $t_1>0$ be small and, for $\mu\in(0,1]$, let us define
$\boQ_{par,\mu}^n=K_n\times[0,\mu t_1]$. $t_1$ is chosen such that
$\boQ_{par,1}^n\subset \boD_n$ and $\boQ_{par,1}^n\cap
\boQ_{uns}^{n,1}=\emptyset$. $\boQ_{par,1}^n$ is foliated by the
equidistant surfaces to $\Sigma_0$ and we have $\Div \xi\le 2H_0$ since the
$H\le H_0$ hypothesis is satisfied.

\subsubsection{Construction of compact stable constant mean curvature
surfaces}

With the notations of the preceding subsection, we have the following lemma.

\begin{lem}\label{minimizing}
There exists $\eps_2\in(0,\eps_1)$ and $p_0\in S$ such that, for large $n$,
there exists a stable constant mean curvature $H_0$ embedded surface $S_n'$
in $(\boD_{n+1}\cap K_n\times\R)\setminus (\boK_{bar}\cup
\boQ_{par,1/2}^{n+1})$ with boundary in
$\partial K_n\times\R$ and $K_n\times\{\eps_2\}$ and $S_n'\cap
[\pi(p_0),p_0]\neq\emptyset$. Moreover the surfaces $S'_n$ are
well oriented \textit{i.e.} the mean curvature vector points into the
connected component of $(\boD_{n+1}\cap K_n\times\R)\setminus S'_n$ which
contains $K_n\times{0}$ and the surfaces $S'_n$, $n$ large, satisfy a
uniform local area estimate of at most one leaf.
\end{lem}

Before the proof of Lemma~\ref{minimizing}, let us explain why we
introduced the subsets $\boK_{bar}$, $\boQ_{bar,t}$,
$\boQ_{uns}^{n,1}$, $\boQ_{x,\delta}$ and $\boQ_{par,\mu}^n$. In fact the
subsets
$\boQ_{bar,t}$, $\boQ_{uns}^{n,1}$, $\boQ_{x,\delta}$ and
$\boQ_{par,\mu}^n$ are used as barriers to prevent the surface $S_n'$ from
touching $\boK_{bar}$,
$\boS_n$ and $\Sigma_0$. So $\boQ_{uns}^{n,1}$ and $\boQ_{par,\mu}^n$
are used to prescribe the boundary of $S_n'$. Once we have the sequence
$S_n'$, we construct $S'$ as the limit of this sequence. We then use
$\boK_{bar}$ as a barrier to control the possible limits of the sequence.

Let us come back to the proof.

\begin{proof}[Proof of Lemma~\ref{minimizing}]

Let $\boF$ be the family of open domains $\boQ$ in
$\boD_{n+1}\setminus \boK_{bar}$ with rectifiable boundary such that
$\boS_{n+1}\subset \partial \boQ$. In the following, $\boS$ will denote the
complement of $\boS_{n+1}$ in $\partial \boQ$.
On $\boF$, we define the functional:
$$
F(\boQ)=A(\partial\boQ) +2H_0 V(\boQ)
$$
where $V(\boQ)$ is the volume of $\boQ$ and $A(\partial\boQ)$ is the
$\boH^2$ measure of $\partial\boQ$. We recall that $A(\partial\boQ)$
is also the mass of the current $\partial[\boQ]$, it is interpreted as
the area of the boundary of $\boQ$. The idea is to find $\boQ_0\in\boF$
which minimizes $F$ in $\boF$ then the part of the boundary of $\boQ_0$ in
$\boD_{n+1}$ will be the surface $S'_n$ we look for.

\begin{claim}\label{claim1}
Let $\boQ$ be in $\boF$.
\begin{enumerate}
\item If $\boQ\cap\boQ_{bar,2\eta/3}\neq \emptyset$, there
exists $t\in[2\eta/3,\eta]$ such that
$\boQ\setminus\boQ_{bar,t}\in\boF$ and $F(\boQ\setminus
\boQ_{bar,t})\le F(\boQ)$.
\item If $\boQ\cap\boQ_{par,1/2}^{n+1}\neq \emptyset$,
there exists $\mu\in[1/2,1]$ such that $\boQ\setminus
\boQ_{par,\mu}^{n+1}\in\boF$ and $F(\boQ\setminus
\boQ_{par,\mu}^{n+1})\le F(\boQ)$.
\end{enumerate}
\end{claim}
\begin{proof}[Proof of Claim~\ref{claim1}]
Let $\boQ$ be in $\boF$ and assume that $\boQ\cap\boQ_{bar,2\eta/3}\neq
\emptyset$ as in Assertion \textit{1}. 

Since $\partial\boQ$ has finite $\boH^2$ measure, the coarea formula
implies that there exists $t\in[2\eta/3,\eta]$ such that
$\boH^1(\boS\cap\boC_{\eta,t})<+\infty$. Thus $\boH^2(\boS\cap
\boC_{\eta,t})=0$: this set is negligible in the following
computations.

First $\boQ\cap\boQ_{bar,t}\neq \emptyset$ has a rectifiable boundary, thus
applying Equation~\eqref{forstokes} of Subsection~\ref{stokes} with $\Div
\xi_{bar}\le 2H_0$, we have:
\begin{align*}
2H_0V(\boQ\cap\boQ_{bar,t})&\ge \int_{\boQ\cap\boQ_{bar,t}}\Div \xi_{bar}
\dd\boL_{\dd s^2}\\
&\ge \int_{\partial(\boQ\cap\boQ_{bar,t})}\la
\xi_{bar}(x),n(\boQ\cap\boQ_{bar,t},x)\ra \dd\boH^2_{\dd s^2}\\
&\ge \int_{\boQ\cap\boC_{\eta,t}}\la
\xi_{bar}(x),n(\boQ\cap\boQ_{bar,t},x)\ra
\dd\boH^2_{\dd s^2} +\int_{\overset{\circ}{\boQ_{bar,t}}\cap\boS}\la
\xi_{bar},n(\boQ\cap\boQ_{bar,t},x)\ra \dd\boH^2_{\dd s^2}.
\end{align*}
We notice that the computation are made with respect to the metric $\dd
s^2$ and results of Subsection~\ref{stokes} are still valid in this
setting.
On $\boC_{\eta,t}\cap \boQ$, we have $\xi_{bar}(x)=n(\boQ\cap\boQ_{bar},x)$
everywhere, thus:
\begin{align*}
A(\boQ\cap\boC_{\eta,t})&=\int_{\boQ\cap\boC_{\eta,t}} \la
\xi_{bar},n(\boQ\cap\boQ_{bar},x)\ra \dd\boH^2_{\dd s^2}\\
&\le -\int_{\overset{\circ}{\boQ_{bar,t}}\cap\boS}\la
\xi_{bar},n(\boQ\cap\boQ_{bar},x)\ra \dd\boH^2_{\dd s^2}
+2H_0V(\boQ\cap\boQ_{bar})\\
&\le A(\overset{\circ}{\boQ_{bar,t}} \cap\boS)+2H_0V(\boQ\cap\boQ_{bar}).
\end{align*}
This implies that
\begin{align*}
F(\boQ\setminus\boQ_{bar,t})&=A(\partial\boQ)-A(\overset{\circ}{\boQ_{
par , t }} \cap\boS)+ A(\boQ\cap\boC_{\eta,t})+
2H_0(V(\boQ)-V(\boQ\cap\boQ_{bar}))\\
&\le F(\boQ).
\end{align*}
Assertion \textit{1} is then proved. Assertion \textit{2} follows from
the same arguments.
\end{proof}

Let $K_{n+1/2}$ be a compact subset of $\Sigma$ such that
$$
K_n\subset \overset{\circ}{K}_{n+1/2}\subset K_{n+1/2}\subset
\overset{\circ}{K}_{n+1}.
$$

\begin{claim}\label{claim2}
Let $\boQ$ be in $\boF$. If $\boQ_{uns}^{n+1,1}\cap
(K_{n+1/2}\times[0,\eps_1]) \not\subset \boQ$, there exists $\boQ'\in\boF$
such that $\boQ_{uns}^{n+1,1}\cap (K_{n+1/2}\times [0,\eps_1]) \subset
\boQ'$ and $F(\boQ')\le F(\boQ)$.
\end{claim}
\begin{proof}[Proof of Claim~\ref{claim2}]
Let $\boQ$ be in $\boF$ as in the claim. The subset $\boQ_{uns}^{n+1,1}\cap
(K_{n+1/2}\times [0,\eps_1])$ is compact so there exists a finite
number of points $x_i \in \boS_{t_0}^{n+1}$ such that 
$$
\boQ_{uns}^{n+1,1}\cap (K_{n+1/2}\times [0,\eps_1]) \subset
\bigcup_i\boQ_{x_i,3/2}
$$

As in proof of Claim~\ref{claim1}, there is $\delta_1\in[3/2,2]$ such that
$\boH^2(\boS\cap \boS_{\delta_1 t_{x_1}}^{x_1})=0$. We denote
$\boO_1=(\boQ_{x_1,\delta_1}\cap \boD_{n+1})\setminus \boQ$. The
boundary of $\boO_1$ is composed of a part
$\partial_1\boO_1=\boS\cap \overset{\circ}{\boQ_{x_1,\delta_1}}$, a second 
part $\partial_2\boO_1\subset \boS_{\delta_1t_{x_1}}^{x_1}$ in the
complement of $\overline{\boQ}$ and a third one of vanishing $\boH^2$
measure. In $\boQ_{x_1,\delta_1}$, we have the unit vector field
$\xi_{x_1}$ which satisfies $\Div \xi_{x_1}<-2H_0$. Then:
\begin{align*}
2H_0V(\boO_1)&\le -\int_{\boO_1} \Div \xi_{x_1}\\
&\le -\int_{\partial\boO_1}\la \xi_{x_1}(x),n(\boO_1,x)\ra\\
&\le -\int_{\partial_2\boO_1}\la \xi_{x_1}(x),n(\boO_1,x)\ra
-\int_{\partial_1\boO_1}\la \xi_{x_1}(x),n(\boO_1,x)\ra.
\end{align*}
where, for simplicity, we have omitted to write the measures. On
$\partial_2\boO_1$,
$\xi_{x_1}=n(\boO_1,x)$ thus
\begin{align*}
2H_0V(\boO_1)
+A(\partial_2\boO_1)&=2H_0V(\boO_1)+\int_{\partial_2\boO_1}\la
\xi_{x_1}(x),n(\boO_1,x)\ra\\
&\le -\int_{\partial_1\boO_1}\la \xi_{x_1}(x),n(\boO_1,x)\ra\\
&\le A(\partial_1\boO_1)
\end{align*}
The interior $\boQ_1$ of $\boQ\cup\boO_1$ is an element of $\boF$ (the
boundary is still rectifiable) and
\begin{align*}
F(\boQ_1)&=2H_0(V(\boQ)+V(\boO_1)+A(\partial\boQ)+A(\partial_2\boO_1)
-A(\partial_1\boO_1)\\
&\le F(\boQ)
\end{align*}
Now considering $\boO_2=(\boQ_{x_2,\delta_2}\cap \boD_{n+1})\setminus
\boQ_1$ and $\boQ_2$ the interior of $\boQ_1\cup \boO_2$, we prove by the
same argument that $\boQ_2\in\boF$ and $F(\boQ_2)\le F(\boQ_1)$. Doing this
a finite number of times, we construct the subset $\boQ'$.
\end{proof}

Let us now consider $(\boQ_k)_k$ a minimizing sequence for $F$. Because of
the claims, we can assume that the sequence satisfies $\boQ_k\cap
\boQ_{bar,2\eta/3}=\emptyset$, $\boQ_k\cap
\boQ_{par,1/2}^{n+1}=\emptyset$ and
$\boQ_{uns}^{n+1,1} \cap (K_{n+1/2}\times [0,\eps_1])\subset \boQ_k$. By
the compactness theorem for integral
currents (see Theorem 5.5 in \cite{Mor}), there is
$\boQ_\infty$ a cluster point of the sequence for the
flat topology. As a limit of a subsequence of $(\boQ_k)_k$, $\boQ_\infty$
is a domain in $\boD_{n+1}$ with a rectifiable boundary such that
$\boQ_\infty\cap \boQ_{bar,2\eta/3}=\emptyset$, $\boQ_\infty\cap
\boQ_{par,1/2}^{n+1}=\emptyset$ and $\boQ_{uns}^{n+1,1}\cap (K_{n+1/2}
\times
[0,\eps_1]) \subset \boQ_\infty$. Moreover $\boQ_\infty$ minimizes $F$
since the area functional $A(\partial\boQ)$ is lower semi-continuous for
the flat convergence and $V(\boQ)$ is the integral over
$\boQ$ of the volume differential form. Since $\boQ_\infty$ minimizes $F$,
the part of $\partial\boQ_\infty$ inside the interior of $\boD_{n+1}$ is a
local isoperimetric surface in the sense of \cite{Mor1}, by regularity
theory (see Corollary 3.7 in \cite{Mor1}) we obtain that this part of
$\partial\boQ_\infty$ is a smooth surface which we denote by $\boS_{n+1}$.
Since $\boQ_\infty$ minimizes $F$, $\boS_{n+1}$ has constant mean curvature
$H_0$  with mean curvature vector pointing outside of $\boQ_\infty$ and it
is stable (see computations in \cite{BaDoCaEs}). Since $\boQ_\infty\cap
\boQ_{bar,2\eta/3}=\emptyset$, $\boQ_\infty\cap
\boQ_{par,1/2}^{n+1}=\emptyset$ and $\boQ_{uns}^{n+1,1} \cap
(K_{n+1/2}\times
[0,\eps_1])\subset \boQ_\infty$, the part of the boundary of
$\boS_{n+1}$ in $K_{n+1/2}\times [0,\eps_1]$ is only in
$K_{n+1/2}\times\{\eps_1\}$ (here we speak about
a non necessarily regular boundary).

Once all the surfaces $\boS_{n+1}$ are constructed, we choose
$\eps_2<\eps_1$ a regular value of the distance function for all the
$\boS_{n+1}$ and we define $S_n'=\boS_{n+1}\cap
(K_n\times[0,\eps_2])$. We notice that $S_n'$ may be empty for small $n$ if
$\eps_2$ is too small; but, for large $n$, $S_n'\neq \emptyset$. Let
$p_0\in S\cap M_+(\eps_2)$ be a point such the geodesic arc
$[p_0,\pi(p_0)]$ does not meet the surface $S$. For $n$ large enough
$\pi(p_0)\in\overline{\boD_n}$ and $p_0\not \in \boD_n$, this implies that
$S_n'\cap [p_0,\pi(p_0)]\neq\emptyset$. These surfaces $S_n'$ are in fact
the ones we want to construct. First the surface is
well oriented since it is a part of the boundary of $\boQ_\infty$. For the
area estimate, let us consider a point $p$ in $D$ and $P$ a plane in the
tangent space. Since $\cup_n\boD_n=D$ and the sequence $(\boD_n)_n$ is
increasing with respect to compact subsets, there is $t_0$ and $n_0$
such that, for $t\le t_0$ and $n\ge n_0$, $O_{p,P}(t)$ is a subset of
$\boD_n$. Since $\boQ_\infty$ minimizes $F$ we have
$F(\boQ_\infty)\le F(\boQ_\infty\setminus O_{p,P}(t))$ and
$F(\boQ_\infty)\le F(\boQ_\infty\cup O_{p,P}(t))$, this implies
that:
\begin{gather*}
A(S'_n\cap O_{p,P}(t))+2H_0 V(\boQ_\infty\cap O_{p,P}(t))\le A(\partial
O_{p,P}(t)\cap \boQ_\infty)\\
A(S'_n\cap O_{p,P}(t))\le A(\partial O_{p,P}(t)\setminus \boQ_\infty)
+2H_0 V(O_{p,P}(t)\setminus \boQ_\infty)
\end{gather*}
Thus, taking the sum and dividing by two, 
$$
A(S'_n\cap O_{p,P}(t))\le A(\partial O_{p,P}(t))/2 +H_0 V(O_{p,P}(t)) 
$$
which is uniformly less that $2\beta\pi t^2$ for some $\beta<1$ and $t$
small because of the asymptotic behaviour  of $A(\partial O_{p,P}(t))$
and $V(O_{p,P}(t))$ (see Subsection~\ref{areabound}).
\end{proof}

\subsubsection{Construction of the surface $S'$}

The last step of the proof of Theorem~\ref{constab} is to obtain a limit
to the sequence $(S'_n)_n$. We choose $\eps_3$ less than $\eps_2$ and
we consider $k\in\N$. For every $n\ge k+1$ and $p\in S'_n\cap
(K_k\times[0,\eps_3])$ the distance from $p$ to the boundary of $S_n$
is
bounded from below by a constant depending only on $k$ and $\eps_3$.
From the stability of $S_n$, this implies that the norm of second
fundamental form of $S'_n$ is bounded in $K_k\times[0,\eps_3]$. Besides
the sequence $(S'_n)_n$ satisfies to a uniform local area estimate. The
curvature and the area estimate implies that the sequence of surface has a
subsequence that converge to a stable cmc $H_0$ surface in
$K_k\times[0,\eps_3]$. Because of the area estimate, the convergence has
multiplicity one and the limit surface is embedded. Since the surfaces
$S'_n$ cuts the geodesic arc $[\pi(p_0),p_0]$ we can assume that this is
also
the case for this limit surface. Then by a diagonal process, we obtain a
stable cmc $H_0$ surface $S_\infty$ in $\Sigma\times[0,\eps_3]$. We have
$S_\infty\cap [\pi(p_0),p_0]$ thus $S_\infty\not \subset \Sigma_{\eps_3}$.
Moreover $S_\infty$ is well oriented as limit of well oriented surfaces.

One thing we have to check is that $S_\infty$ is in fact in $\Sigma\times
(0,\eps_3]$. If it is not the case, $S_\infty$ touches $\Sigma_0$ and by
the maximum principle we have $S_\infty=\Sigma_0$. By construction, the
sequence $S'_n$ never enters in $\boK_{bar}$ so it is the same for
$S_\infty$ and we obtain $S_\infty\neq \Sigma_0$. 

Moreover $S_\infty$ is not included in an equidistant surface
$\Sigma_t$.
By construction, $S_\infty$ is between $\Sigma_0$ and $S$ and $\inf_S
\d=0$, this implies that $\d$ can not be constant along $S_\infty$. 

Now choosing $\eps'$ a regular value of the distance function $\d$ on
$S_\infty$ (we assume that $\eps'$ is part of the image of $\d$ along
$S_\infty$), we can consider $S'= S_\infty \cap \Sigma\times [0,\eps']$:
$S'$ then has its non empty boundary in $\Sigma\times\{\eps'\}$. $S'$ is
then a complete stable cmc $H_0$ surface which is properly embedded in
$\Sigma\times[0,\eps']$. Moreover $S'$ is well oriented and $\d$ is not
constant along $S'$.

\subsection{$H_0=0$}

In this case, the cases \textit{1} and \textit{2} of
Theorem~\ref{constab} are the same, so assume that $S\looparrowright
M_+(\eps)$. The proof is essentially the same, the only difference comes
from the fact that the ``well oriented" hypothesis has no more meaning.

So as above we define, $\boK_{bar}$, $\boQ_{bar,t}$ and $\xi_{bar}$ such
that $\Div\xi_{bar} \le 0$. This gives a $\eps_1$.

We introduce the compact $K_n$ and the domain $D_n$. As above we assume
the instability of $S_n$ and consider $\phi$ such that
$L_\phi=\lambda_1\phi$, $v$ such that $-Lv=1$ and $u=\phi+av$. Let
$N(x)$ be the unit normal to $S_n$. For $t_0>0$ we define, on
$S_n\times[-t_0,t_0]$, the map $F(x,t)=\exp_x(tu(x)N(x))$ and assume
that $t_0$ is small enough to ensure that $F$ is an immersion.
$S_n\times[-t_0,t_0]$ with the metric $F^*\dd s^2$ is foliated by
$\boS_t^n=S_n\times\{t\}$. Because of $-Lu=-\lambda_1\phi+a>0$, if $t_0$
is chosen small enough, the mean curvature vector of $F(\boS_{t_0}^n)$
and $F(\boS_{-t_0}^n)$ is non vanishing and points ``outside"
$F(S_n\times[-t_0,t_0])$.

Thus for any $x\in \boS_{t_0}^n\cup \boS_{-t_0}^n$ we can define as above
$\boQ_{x,\delta}$ and $\xi_x$ with $\Div \xi_x <0$. Then we
define $\boQ_{uns}^{n,0}=F(S_n\times[-t_0,t_0])\cap D_n$,
$\boD_n=D_n\setminus \boQ_{uns}^{n,0}$ and
$$
\boQ_{uns}^{n,1}=\boD_n\cap \bigcup_{x\in\boS_{t_0}^n\cup
\boS_{-t_0}^n}\boQ_{x,1}
$$

With these notations, the end of the proof is the same.

\subsection{$S\looparrowright M_-(\eps)$ and $H_0>0$}

When $S\looparrowright M_-(\eps)$, the differences comes from the fact
that the surface is not assumed to be well oriented.

As above, we define, $\boK_{bar}$, $\boQ_{bar,t}$ and $\xi_{bar}$ such that
$\Div\xi_{bar} \le -2H_0$. This gives $\eps_1$. We introduce the compact
subsets $K_n$ and the domain $D_n$.

We use the instability of $S$ to define $\phi$ such that
$L_\phi=\lambda_1\phi$, $v$ such that $-Lv=1$ and $u=\phi+av$. If $N(x)$
is the unit normal to $S_n$ such that the mean curvature vector is $2H_0
N(x)$, we define $F(x,t)=\exp_x(tu(x)N(x))$ on $S_n\times [-t_0,0]$ with
$t_0>$ small so that $F$ is an immersion. $S_n\times[-t_0,0]$
with the metric $F^*\dd s^2$ is foliated by
$\boS_t^n=S_n\times\{t\}$ and we extend to $S_n\times[-t_0,0]$ the
definition of $N$ as the unit normal vectorfield to the surfaces
$\boS_t^n$.

Since $-Lu=-\lambda_1\phi+a>0$, if $t_0$ is small enough, the mean
curvature of $F(\boS_{-t_0}^n)$ computed with respect to $N$ is less than
$H_0$. We define $\boQ_{uns}^{n,0}=F(S_n\times[-t_0,0]) \cap D_n$ and
$\boD_n=D_n\setminus \boQ_{uns}^{n,0}$.

As above for any $x\in \boS_{-t_0}^n$, we consider $D_x\subset
\boS_{-t_0}^n$ a small open geodesic disk which is embedded in
$M_+(\eps_1)$ by $F$. Let $\psi$ be a smooth function on $\boS_{t_0}^n$
vanishing outside $D_x$ and positive in $D_x$. We then define on
$\boS_{-t_0}^n\times[0,2t_x]$ 
$$
G(p,t)=\exp_{F(p)}(-t\psi(p)F_*(N(p)))
$$
If we choose $t_x$ small, we can assume that $G$ is an embedding on
$D_x\times[0,2t_x]$. In $G(D_x\times[0,2t_x])$, we define $\xi_x$ the unit
vector field normal to the embedded surfaces $\boS_t^x=G(D_x\times\{t\})$
with $\xi_x=-F_*(N)$ along $\boS_0^x$. Since the mean curvature of
$\boS_0^x$ is less than $2H_0$, if $t_x$ is small, we can
assume that $\Div \xi_x<2H_0$.

We denote $\boQ_{x,\delta}=G(D_x\times[0,\delta t_x])$ for $\delta\in
[1,2]$. We also define $$
\boQ_{uns}^{n,1}=\boD_n\cap\bigcup_{x\in
\boS_{t_0}^n}\boQ_{x,1}.
$$

Since the surface $S_n$ can be not well oriented, we need to introduce
new
barriers. For any $x\in S_n$, we consider $\Delta_x\subset
S_n$ a small open geodesic disk which is embedded in
$M_+(\eps_1)$. Let $\eta$ be a smooth function on $S_n$ vanishing outside
$\Delta_x$ and positive in $\Delta_x$. We define on $S_n\times[0,2t_x]$
$$
H(p,t)=\exp_p(t\eta(p)N(p))
$$
with $t_x$ small enough such that $H$ is an embedding on
$\Delta_x\times[0,2t_x]$. In $H(\Delta_x\times[0,2t_x])$, we define
$\xi_{ori}^x$ the unit vector field normal to the surfaces
$\boS_{ori,t}^x=H(\Delta_x\times\{t\})$ with $\xi_{ori}^x=N$ along
$\Delta_x$. Since the mean curvature vector of $S_n$ is $H_0N$, if $t_x$
is small enough we have $\Div \xi_{ori}^x<2H_0$.

We denote $\boQ_{ori}^{x,\nu}=H(\Delta_x\times[0,\nu t_x])$ and
$\overset{\circ}{\boQ_{ori}^{x,\nu}}=H(\Delta_x\times[0,\nu t_x))$,
 for $\nu\in[1,2]$ and :
$$
\boQ_{ori}^{n,1}=\boD_n\cap\bigcup_{x\in
S_n}\boQ_{ori}^{x,1}
$$

Let $\boS_n$ be part of the boundary of $\boD_n$ in $F(\boS_{-t_0}^n)\cup
S_n$; this the part of $\partial\boD_n$ not in $\Sigma_0$ and $\partial
K_n\times\R$. Any point in $\boS_n$ is at positive distance  from
$\boD_n\setminus (\boQ_{uns}^{n,1}\cup \boQ_{ori}^{n,1})$.

Then we define $\boQ_{par,\mu}^n=K_n\times[-\mu t_1,0]$ and introduce
$\xi_{par}=-\xi$. We have $\Div \xi_{par}=-\Div\xi \le -2H_0$ because
of the $H\ge H_0$ hypothesis in $M_-(\eps)$.

An equivalent of Lemma~\ref{minimizing} can be proved. The idea is now to
minimize the functional $F(\boQ)=A(\partial\boQ)-2H_0V(\boQ)$ where
$\boQ\in\boF$
and $\boF$ is the same family of domains in $\boD_{n+1}$. We first
remark that Claim~\ref{claim1} is still true. Claim~\ref{claim2} is
replaced by
\begin{claim}\label{claim3}
Let $\boQ$ be in $\boF$.
\begin{enumerate}
\item If $\boQ_{uns}^{n+1,1}\cap
(K_{n+1/2}\times[0,\eps_1]) \not\subset \boQ$, there exists
$\boQ'\in\boF$ such that $\boQ_{uns}^{n+1,1}\cap (K_{n+1/2}\times
[0,\eps_1]) \subset \boQ'$ and $F(\boQ')\le F(\boQ)$.
\item If $\boQ_{ori}^{n+1,1}\cap (K_{n+1/2}\times[0,\eps_1]) \not\subset
\boQ$, there exists $\boQ'\in\boF$ such that $\boQ_{ori}^{n+1,1}\cap
(K_{n+1/2}\times [0,\eps_1]) \subset \boQ'$ and $F(\boQ')\le F(\boQ)$.
\end{enumerate}
\end{claim}

\begin{proof}[Proof of Claim~\ref{claim3}]
The proof of both items are the same so let us prove the second one. Let
$\boQ$ be in $\boF$ and $\boQ_{ori}^{n+1,1}\cap (K_{n+1/2}\times[0,\eps_1])
\not\subset \boQ$. As in the proof of Claim~\ref{claim2}, the subset
$\boQ_{ori}^{n+1,1}\cap (K_{n+1/2}\times [0,\eps_1])$ is compact so there
exists a finite number of points $x_i \in S_{n+1}$ such that 
$$
\boQ_{ori}^{n+1,1}\cap (K_{n+1/2}\times [0,\eps_1]) \subset
\bigcup_i\boQ_{ori}^{x_i,3/2}.
$$

As in claims~\ref{claim1} and \ref{claim2}, there is $\nu_1\in[3/2,2]$ such
that $\boH^2(\boS\cap \boS_{ori,\nu_1t_{x_1}}^{x_1}) =0$. Then we denote
$\boO_1=(\boQ_{ori}^{x_1,\nu_1}\cap \boD_{n+1})\setminus \boQ$. The
boundary of $\boO_1$ is composed of a part $\partial_1\boO_1=
\boS\cap\overset{\circ}{\boQ_{ori}^{x_1,\nu_1}}$ and a second part
$\partial_2\boO_1\subset \boS_{ori,\nu_1 t_{x_1}}^{x_1}$ in the complement
of $\overline{\boQ}$ and a third part of vanishing $\boH^2$ measure. In
$\boQ_{ori}^{x_1,\nu_1}$, we have the unit vector field $\xi_{ori}^{x_1}$
which satisfies $\Div \xi_{ori}^{x_1}<2H_0$. Then
:
\begin{align*}
2H_0V(\boO_1)&\ge \int_{\boO_1} \Div \xi_{ori}^{x_1}\\
&\ge \int_{\partial\boO_1}\la \xi_{ori}^{x_1},n(\boO_1,x\ra\\
&\ge \int_{\partial_2\boO_1}\la \xi_{ori}^{x_1},n(\boO_1,x)\ra
+\int_{\partial_1\boO_1}\la \xi_{ori}^{x_1},n(\boO_1,x)\ra .
\end{align*}
On $\partial_2\boO_1$, $\xi_{ori}^{x_1}=n(\boO_1,x)$ thus
\begin{align*}
A(\partial_2\boO_1)-2H_0V(\boO_1)&=\int_{\partial_2\boO_1} \la
\xi_{ori}^{x_1},n(\boO_1,x)\ra-2H_0V(\boO_1)\\
&\le -\int_{\partial_1\boO_1}\la \xi_{ori}^{x_1},n(\boO_1,x)\ra\\\
&\le A(\partial_1\boO_1).
\end{align*}
$\boQ_1=\boQ\cup\boO_1$ is an element of $\boF$  since the boundary is
still rectifiable and
\begin{align*}
F(\boQ\cup\boO_1)&=-2H_0(V(\boQ)+V(\boO_1))+A(\partial\boQ)+A(\partial_2
\boO_1)-A(\partial_1\boO_1)\\
&\le F(\boQ).
\end{align*}
Repeating this a finite number of times we construct the subset $\boQ'$.
\end{proof}

As in proof of Lemma~\ref{minimizing}, we obtain a minimizer $\boQ_\infty$
and a smooth surface $\boS_{n+1}$ which gives us $S'_n$. The uniform
area estimate is also proved by the same way.

Once the sequence $S'_n$ is constructed, the end of the proof of
Theorem~\ref{constab} is the same as in the first case.

%% file: fighalfspace1.pstex_t
\begin{picture}(0,0)%
\includegraphics{fighalfspace1.pstex}%
\end{picture}%
\setlength{\unitlength}{4144sp}%
\begingroup\makeatletter\ifx\SetFigFontNFSS\undefined%
\gdef\SetFigFontNFSS#1#2#3#4#5{%
  \reset@font\fontsize{#1}{#2pt}%
  \fontfamily{#3}\fontseries{#4}\fontshape{#5}%
  \selectfont}%
\fi\endgroup%
\begin{picture}(4932,1447)(706,-3023)
\put(721,-2131){\makebox(0,0)[lb]{\smash{{\SetFigFontNFSS{14}{16.8}{\rmdefault}{\mddefault}{\updefault}{\color[rgb]{0,0,0}$S$}%
}}}}
\put(3736,-1771){\makebox(0,0)[lb]{\smash{{\SetFigFontNFSS{14}{16.8}{\rmdefault}{\mddefault}{\updefault}{\color[rgb]{0,0,0}$\boQ_{bar,\eta}$}%
}}}}
\put(4366,-2941){\makebox(0,0)[lb]{\smash{{\SetFigFontNFSS{14}{16.8}{\rmdefault}{\mddefault}{\updefault}{\color[rgb]{0,0,0}$\boC_{\eta,t}$}%
}}}}
\put(4321,-2131){\makebox(0,0)[lb]{\smash{{\SetFigFontNFSS{14}{16.8}{\rmdefault}{\mddefault}{\updefault}{\color[rgb]{0,0,0}$\boK_{bar}$}%
}}}}
\end{picture}%

%% file: fighalfspace2.pstex_t
\begin{picture}(0,0)%
\includegraphics{fighalfspace2.pstex}%
\end{picture}%
\setlength{\unitlength}{4144sp}%
\begingroup\makeatletter\ifx\SetFigFontNFSS\undefined%
\gdef\SetFigFontNFSS#1#2#3#4#5{%
  \reset@font\fontsize{#1}{#2pt}%
  \fontfamily{#3}\fontseries{#4}\fontshape{#5}%
  \selectfont}%
\fi\endgroup%
\begin{picture}(7449,2949)(439,-3448)
\put(3218,-2926){\makebox(0,0)[lb]{\smash{{\SetFigFontNFSS{14}{16.8}{\rmdefault}{\mddefault}{\updefault}{\color[rgb]{0,0,0}$\boD_n$}%
}}}}
\put(6991,-1471){\makebox(0,0)[lb]{\smash{{\SetFigFontNFSS{14}{16.8}{\rmdefault}{\mddefault}{\updefault}{\color[rgb]{0,0,0}$S_n$}%
}}}}
\put(7036,-2206){\makebox(0,0)[lb]{\smash{{\SetFigFontNFSS{14}{16.8}{\rmdefault}{\mddefault}{\updefault}{\color[rgb]{0,0,0}$\boS_{t_0}^n$}%
}}}}
\put(608,-2356){\makebox(0,0)[lb]{\smash{{\SetFigFontNFSS{14}{16.8}{\rmdefault}{\mddefault}{\updefault}{\color[rgb]{0,0,0}$\boQ_{uns}^{n,0}$}%
}}}}
\put(4764,-2596){\makebox(0,0)[lb]{\smash{{\SetFigFontNFSS{14}{16.8}{\rmdefault}{\mddefault}{\updefault}{\color[rgb]{0,0,0}$\boQ_{x,2}$}%
}}}}
\put(4742,-1777){\makebox(0,0)[lb]{\smash{{\SetFigFontNFSS{14}{16.8}{\rmdefault}{\mddefault}{\updefault}{\color[rgb]{0,0,0}$\boQ_{x,1}$}%
}}}}
\end{picture}%

%% file: halfspace5.tex
In this section we prove our main theorem.

\subsection{Some preliminary computations}

We begin by some computations. Let $\Sigma$ be a Riemannian surface and
$M_\pm(\eps)$ be an
$\eps$-half neighborhood of $\Sigma$. Let $S$ be a constant mean curvature
$H_0$ surface in $M_\pm(\eps)$. We denote by $\nabla$ the connection on
$M_\pm(\eps)$ and we denote by $\widetilde{\nabla}$ and
$\widetilde{\Delta}$ the connection on $S$ and its associated Laplace
operator.

Let $f$ be a function on $\R$, we want to compute $\widetilde{\Delta}
f(\d)$. Along $S$, we denote by $(e_1,e_2,e_3)$ an orthonormal basis of
$TM_\pm(\eps)$ such that $e_3$ is normal to $S$ and the mean curvature
vector to  $S$ is $H_0 e_3$. For any function $g$ defined in $M_\pm(\eps)$
we have :
$$
\widetilde{\Delta}g=\sum_{i=1}^2\la\nabla_{e_i}\nabla g,e_i\ra
+2\la\nabla g,H_0 e_3\ra
$$
Thus if $g=f\circ \d$, we get:
$$
\widetilde{\Delta}f\circ \d=f''(\d)\sum_{i=1}^2\la \nabla
\d,e_i\ra^2+f'(\d)\left( \sum_{i=1}^2\la\nabla_{e_i}\nabla \d,e_i\ra+ 2\la
\nabla \d,H_0 e_3\ra\right)
$$

For any point in $M_\pm(\eps)$, we denote by $(a_1,a_2)$ an orthonormal
basis of $T\Sigma_t$ which diagonalized the shape operator of $\Sigma_t$.
Let $\kappa_1$ and $\kappa_2$ the associated principal curvature such that
$ \nabla_{a_i}\xi=\kappa_i a_i$ for $i=1,2$. $(a_1,a_2,\xi)$ is then an
orthonormal basis of $TM_\pm(\eps)$, we write $a_3=\xi$ and we have
$\nabla_{a_3}a_3=0$. Moreover, we define $(\lambda_i^j)_{1\le i,j\le 3}$
such that 
$$
e_i=\sum_{j=1}^3\lambda_i^j a_j
$$
Using these expressions and working in $M_+(\eps)$ where $\nabla
\d=a_3$, we have:
\begin{align*}
\widetilde{\Delta}f\circ \d&= f''(\d)
({\lambda_1^3}^2+{\lambda_2^3}^2)+f'(\d)\left(\sum_{i=1} ^2\la
\lambda_i^1\kappa_1a_1+\lambda_i^2\kappa_2 a_2,\lambda_i^1 a_1+\lambda_i^2
a_2+\lambda_i^3a_3\ra+2H_0\lambda_3^3\right)\\
&=f''(\d)(1-{\lambda_3^3}^2)+f'(\d)\left( \kappa_1
({\lambda_1^1}^2+{\lambda_2^1}^2)+
\kappa_2({\lambda_1^2}^2+{\lambda_2^2}^2)+2H_0\lambda_3^3\right)\\
&=f''(\d)(1-{\lambda_3^3}
^2)+f'(\d)\left((\kappa_1+\kappa_2)+2H_0\lambda_3^3-\kappa_1{\lambda_3^1}
^2-\kappa_2 {\lambda_3^2}^2\right)
\end{align*}

Since the vector $(\lambda_3^1,\lambda_3^2,\lambda_3^3)$ has norm $1$,
there exists $(\phi,\theta)\in[0,\pi]\times[0,2\pi]$ such that 
$$
(\lambda_3^1,\lambda_3^2,\lambda_3^3)=(\sin\phi\cos \theta,\sin\phi
\sin\theta, -\cos\phi)
$$
The ``$-$'' sign in the last coordinate is there in order to make $\phi$
close to $0$ in the proof below. Besides, if $M_+(\eps)$ satisfies the
$H\le H_0$ hypothesis and $f$ is an increasing function, we obtain:
$$
\widetilde{\Delta}f\circ \d \le f''(\d)
(1-\cos^2\phi)+ f'(\d)(2H_0(1-\cos\phi)- (\kappa_1\cos^2\theta+
\kappa_2\sin^2 \theta)\sin^2\phi)$$
\begin{equation}\label{surharmonic1}
\widetilde{\Delta}f\circ \d \le f''(\d)\sin^2\phi+ f'(\d)(2H_0(1-\cos\phi)-
(\kappa_1\cos^2\theta+ \kappa_2\sin^2 \theta)\sin^2\phi)
\end{equation}

If we work in $M_-(\eps)$, we have $\nabla \d=-a_3$ thus:
\begin{align*}
\widetilde{\Delta}f\circ \d&= f''(\d)
({\lambda_1^3}^2+{\lambda_2^3}^2)-f'(\d)\left(\sum_{i=1} ^2\la
\lambda_i^1\kappa_1a_1+\lambda_i^2\kappa_2 a_2,\lambda_i^1 a_1+\lambda_i^2
a_2+\lambda_i^3a_3\ra+2H_0\lambda_3^3\right)\\
&=f''(\d)(1-{\lambda_3^3}^2)-f'(\d)\left( \kappa_1
({\lambda_1^1}^2+{\lambda_2^1}^2)+
\kappa_2({\lambda_1^2}^2+{\lambda_2^2}^2)+2H_0\lambda_3^3\right)\\
&=f''(\d)(1-{\lambda_3^3}
^2)-f'(\d)\left((\kappa_1+\kappa_2)+2H_0\lambda_3^3-\kappa_1{\lambda_3^1}
^2-\kappa_2 {\lambda_3^2}^2\right)
\end{align*}

The vector $(\lambda_3^1,\lambda_3^2,\lambda_3^3)$ has still norm $1$, so
there exists $(\phi,\theta)\in[0,\pi]\times[0,2\pi]$ such that 
$$
(\lambda_3^1,\lambda_3^2,\lambda_3^3)=(\sin\phi\cos \theta,\sin\phi
\sin\theta, -\cos\phi)
$$
If $M_-(\eps)$ satisfies the $H\ge H_0$ hypothesis and $f$ is an increasing
function, we obtain:
\begin{align*}
\widetilde{\Delta}f\circ \d &\le f''(\d) (1-\cos^2\phi)+ f'(\d)
(2H_0(\cos\phi-1)+ (\kappa_1\cos^2\theta+ \kappa_2\sin^2
\theta)\sin^2\phi)\\
&\le f''(\d)\sin^2\phi+ f'(\d)(2H_0(\cos\phi-1)+
(\kappa_1\cos^2\theta+ \kappa_2\sin^2 \theta)\sin^2\phi)
\end{align*}
Since $\cos\phi-1\le 0$, we get:
\begin{equation}\label{surharmonic2}
\widetilde{\Delta}f\circ \d \le \left(f''(\d)+f'(\d)(
\kappa_1\cos^2\theta+ \kappa_2\sin^2 \theta)\right)\sin^2\phi
\end{equation}

\subsection{The main theorem}

Let us now state and prove our main result.

\begin{thm}\label{halfspace}
Let $(\Sigma,\dd \sigma_0^2)$ be a complete orientable Riemannian surface,
$\eps$ be  positive and $H_0$ non-negative. Let $M_\pm(\eps)$ be an inside
or outside $\eps$-half neighborhood of $\Sigma$. We consider a
properly immersed constant mean curvature $H_0$ surface $S$ in
$M_\pm(\eps)$ with possibly non-empty boundary in $\Sigma_\eps$ and
$S\subset M_\pm^*(\eps)$. 

We assume that $(\Sigma,\dd\sigma_0^2)$ is parabolic. We also assume
that $M_\pm(\eps)$ is regular.

\begin{enumerate}
\item If $S\looparrowright M_+(\eps)$ is well oriented and $M_+(\eps)$
satisfies the $H\le H_0$ hypothesis. The distance function $\d$ is constant
on $S$.
\item If $S\looparrowright M_-(\eps)$ and $M_-(\eps)$ satisfies the $H\ge
H_0$ hypothesis. The distance function $\d$ is constant on $S$.
\end{enumerate}

\end{thm}

Theorem~\ref{halfspace} says that the equidistant surfaces are the
only possible constant mean curvature $H_0$ surfaces in $M_\pm(\eps)$ (with
good orientation in $M_+(\eps)$). If no equidistant surface has mean
curvature $H_0$, no cmc $H_0$ surface exists in $M_\pm(\eps)$.

As for Theorem~\ref{constab}, the proof of both cases are very similar
so we will mainly focus on the first one.

\subsubsection{$S\looparrowright M_+(\eps)$ and $H_0>0$}

Let us consider $S$ in $M_+(\eps)$ and assume that $S$ is not in one
equidistant surface $\Sigma_t$ ($\d$ is not constant along $S$). Let $\mu<\eps$ be
the lower bound of $\d$ on $S$, we notice that this lower bound
is never reached because of the $H\le H_0$ hypothesis and the maximum
principle.

The space $N_+(\eps-\mu)=\Sigma\times[\mu,\eps]$ with the Riemannian
metric $\dd s^2$ can be seen as an outside $(\eps-\mu)$-half neighborhood
of
$(\Sigma,\dd\sigma_\mu^2)$. Since $M_+(\eps)$ is regular,
$(\Sigma,\dd\sigma_\mu^2)$ is parabolic ($\pi_\mu:\Sigma_\mu
\rightarrow\Sigma_0$ is quasi-isometric) and $N_+(\eps-\mu)$ is regular.

$S$ can be viewed as properly immersed in $N_+(\eps-\mu)$; thus we can
assume that $\inf_S\d=0$ in the statement of Theorem~\ref{halfspace}.

Thus $S\looparrowright M_+(\eps)$ satisfies all the hypotheses of
Theorem~\ref{constab}. So there are $\eps'>0$ and a surface $S'$ properly
immersed in $M_+(\eps')$ with nonempty boundary in $\Sigma_{\eps'}$. $S'$
is
well oriented, has cmc $H_0$ and is stable, moreover the distance function
$\d$ on $S'$ is not constant.

Let $\eps_1$ be less that $\eps'$; for any point in $S'\cap
M_+(\eps_1)$, the geodesic distance to $\partial S'$ is lower bounded by
$\eps'-\eps_1$. Since $S'$ is stable and $M_+(\eps')$ is geometrically
bounded, the norm of the second fundamental form of $S'$ is bounded in
$M_+(\eps_1)$. Choosing $\eps_1$ sufficiently close to $\inf_{S'}\d$, there
is
a constant $c>0$ such that, along $S'\cap M_+(\eps_1)$, $|\la n,\xi\ra|>c$
where $n$ is the normal to $S'$. Thus $\pi$ is a local quasi-isometry from
$S'$ to $\Sigma_0$. 

Let $D$ by the connected component of $M_+(\eps_1)\setminus S'$ which
contains $\Sigma_0$. Let $p$ be in $S'\cap M_+(\eps_1)$ and $q=\pi(p)$,
along the geodesic segment $[q,p]$ there is a point $p'\in S'$ which is the
closest to $q$. $p'$ is in $\partial D$ and we denote by $S''$ the
connected component of $S'\cap M_+(\eps_1)$ which contains $p'$. $S''$ is
not in an equidistant surface to $\Sigma_0$. We notice that since $S'$ is
well
oriented the mean curvature vector at $p'$ points into $D$. Thus $\la
n,\xi\ra\le 0$ at $p'$ which gives $\la n,\xi\ra\le -c$ at $p'$. Since
$S''$ is connected, we get $\la n,\xi\ra\le -c$ along $S''$.

Let us construct on $S''$ a non constant bounded superharmonic function
which does not reach its lower bound on the boundary. Let $K$ be a
real constant and consider the function:
$$
f_K:\R_+\rightarrow \R,\ x\mapsto \frac{1}{K}(1-\exp(-Kx))
$$
We have $f_K'(x)=\exp(-Kx)\ge 0$ so $f_K$ is increasing and
$f_K''(x)+Kf_K'(x)=0$.

Now, we use the computation~\eqref{surharmonic1} with $f=f_K$. On
$S''$ we have $|\la n,\xi\ra|<-c$, this means that $\cos\phi\ge c $ in
\eqref{surharmonic1}. But there exists $A\ge 0$ such that $1-\cos\phi\le
A\sin^2\phi$ when $\cos\phi\ge c$. Then, from \eqref{surharmonic1}, we get:
\begin{align*}
\Delta_{S''}f_K\circ \d&\le f_K''(\d)\sin^2\phi+ f_K'(\d)(2H_0A\sin^2\phi-
(\kappa_1\cos^2\theta+ \kappa_2\sin^2 \theta)\sin^2\phi)\\
&\le \left(f_K''(\d)+ f_K'(\d)(2H_0A-(\kappa_1\cos^2\theta+ \kappa_2\cos^2
\theta))\right)\sin^2\phi
\end{align*}

Since $M_+(\eps)$ is assumed to be regular there is a constant $C$ such
that
$\max(|\kappa_1|,|\kappa_2 |) \le C$. Then considering $K=2H_0A+C$ we get 
\begin{align*}
\Delta_{S''}f_K\circ \d&\le \left(f_K''(\d)+
f_K'(\d)(2H_0A+(C\cos^2\theta+ C\sin^2 \theta))\right)\sin^2\phi\\
&\le \left(f_K''(\d)+ f_K'(\d)(2H_0A+C)\right)\sin^2\phi\\
&\le 0
\end{align*}

$f_K\circ \d$ is then superharmonic on $S''$, bounded since $\d$ is bounded
and $f_K\circ \d\le f_K(\eps_1)=(f_K\circ \d)_{|\partial S''}$. If we
prove that $S''$ is parabolic at infinity we could conclude that
$f_K\circ \d$ is constant and $S''\subset\Sigma_{\eps_1}$; this will give
the contradiction we look for and the first case of Theorem~\ref{halfspace}
will be proved.

First we deal with a special case: $S''$ is embedded. This case is not
necessary for the general one but it explains some ideas. We
have the following claim
\begin{claim}\label{injective}
$\pi$ is injective on $S''$.
\end{claim}
\begin{proof}
Let us assume that there is $p_0$ and $p_1$ in $S''$ such that
$\pi(p_0)=\pi(p_1)$ and $\d(p_0)>\d(p_1)$. Let $\gamma:[0,1]\rightarrow
S''$ be a curve such that $\gamma(0)=p_0$ and $\gamma(1)=p_1$. We denote
$\pi\circ \gamma$ by $\tilde\gamma$. $\tilde\gamma$ is a closed curve in
$\Sigma_0$, so we can extend the definition of $\tilde\gamma$ by
periodicity to $\R_+$. Since $\pi:S''\rightarrow \Sigma_0$ is a local
diffeomorphism, we can extend the definition of $\gamma$ as a lift of
$\tilde\gamma$ to $[0,t_0]$ where $\gamma(t_0)\in\partial S''$ or to
$\R_+$.

We have $\d(\gamma(0))-\d(\gamma(1))>0$ then, for any $t\in[0,t_0-1]$,
$\d(\gamma(t))-\d(\gamma(t+1))>0$ since this quantity never vanishes. Since
$\d(\gamma(t))\le\eps'$, we get $\d(\gamma(t))<\eps'$ for any $t\ge 1$.
Hence $\gamma(t)\notin \partial S''$ for $t\ge 1$ and $\gamma$ is then
defined on $\R_+$. Thus $\gamma(n)$ is a sequence of distinct points in
$S''$ with $\pi(\gamma(n))=\pi(p_0)$. This contradicts the fact that $S''$
is properly embedded and $|\la n,\xi\ra|>c$. The map $\pi$ is then injective on
$S''$. 
\end{proof}

Since $\pi:S''\rightarrow \Sigma_0$ is an injective quasi-isometry and
$\Sigma_0$ is parabolic, $S''$ is parabolic at infinity by
Proposition~\ref{parabolic}; Theorem~\ref{halfspace} would then be
proved.

Let us now write the general case: $S''$ is only immersed.

We recall that $D_0$ is the connected component of $M_+(\eps_1)\setminus
S''$ that contains $\Sigma_0$. The boundary of $D_0$ is composed by
$\Sigma_0$ and a set $S_0$ made of points in $S''$ and $\Sigma_{\eps_1}$
(see Figure~\ref{defS0}).
For any $x$ in $\Sigma_0$, we define $v(x)=\min\{\d(p),p\in
\pi^{-1}(x)\cap (S''\cup\Sigma_{\eps_1})\}$. It is clear that the graph of
$v$, $\{(x,v(x))\in \Sigma\times[0,\eps_1]\}$, is included in $S_0$. In
fact we have equality because of the following claim.

\begin{figure}[h]
\begin{center}
\resizebox{0.8\linewidth}{!}{\input{fighalfspace3.pstex_t}}
\caption{}
\label{defS0}
\end{center}
\end{figure}

\begin{claim}
The function $v$ is continuous.
 \end{claim}
\begin{proof}
If $v$ is not continuous there is a sequence of points $(x_n)$ converging
to $x$ in $\Sigma_0$ such that $\lim v(x_n)=v_0\neq v(x)$. Since
$S''\cup\Sigma_{\eps_1}$ is closed, $(x,v_0)\subset S''\cup
\Sigma_{\eps_1}$ thus $v_0>v(x)$. $(x,v(x))$ is in $S''\cup
\Sigma_{\eps_1}$ thus there is a smooth function $f$ defined in a
neighborhood of $x$ in $\Sigma_0$ such that the graph of $f$ is included
in $S''\cup \Sigma_{\eps_1}$ and $f(x)=v(x)$ (we used the fact that $|\la
n,\xi\ra|>c$ along $S''$). Then $f<v_0$ near $x$ and $v(x_n)\le f(x_n)$
for $n$ large. We get a contradiction.
\end{proof}

In fact near a point $p\in S_0$, $S''$ and $\Sigma_{\eps_1}$ can be viewed
as a finite union of graphs above a small disk in $\Sigma_0$ around
$\pi(p)$. Let us denote the associated functions by $f_i$, then $v=\min_i
f_i$ (in
view of Subsection~\ref{selfinter}, $f_0=\eps_1$ and $f_1,\cdots, f_p$
have constant mean curvature graphs). The
projection map $\pi :S_0\rightarrow \Sigma_0$ is then a homeomorphism.

Let us denote by $O_i$ the connected component of $S_0$ minus the set of
self-intersection points in $S''$ and the set $S''\cap \Sigma_{\eps_1}$
(these are the points where $v$ is given by only one $f_j$).

We denote $\Ome_i=\pi(O_i)\subset\Sigma_0$. By the description made in
Subsection~\ref{selfinter}, the boundary of $O_i$ can be decomposed as
the union
of part $\Gamma_{i,j}$ and a set of vanishing $\boH^1$ measure. The set
$\Gamma_{i,j}$ is the part of $\partial \Ome_i\cap \partial \Ome_j$ where
$\Ome_i$ ``touches''  $\Ome_j$. On $\overline{\Ome_i}$, we
consider the metric $g_i=\pi_*(\dd s^2_{|O_i})$, this metric is well
defined since $\pi$ is smooth on $S''$ and $\Sigma_{\eps_1}$. Moreover
since $\pi$ is quasi-isometric along $S''$ and $\Sigma_{\eps_1}$ there is
$k>0$ such that $\frac{1}{k^2}\dd \sigma_0^2\le g_i\le k^2\dd \sigma_0^2$.

On $\Sigma_0$ we consider the function $u$ defined by $u(p)= f_K\circ
\d(\pi_{|S_0}^{-1}(p))-\eps_1=f_K\circ v-\eps_1$. $u$ is non-positive,
smooth
on each $\Ome_i$ and $\Delta_{g_i}u\le 0$. In fact, in view of its
definition and the definition of $S_0$, $u$ can be interpreted as the
minimum of several superharmonic functions so, in some sense, $u$ is a
superharmonic function. Let us explain how this idea can be used. The
following computation are inspired by \cite{MaPeRo} (see also
\cite{AmCa,BeCaNi}).

Since $(\Sigma_0,\dd\sigma_0^2)$ is parabolic there exists a sequence of
compactly supported smooth functions $(\phi_n)_n$ such that $0\le \phi_n\le
1$, $(\phi_n^{-1}(1))_n$ is a compact exhaustion of $\Sigma_0$ and 
$$
\lim_n\int_{\Sigma_0}\|\nabla_0\phi_n\|_0\dd v_0=0
$$
The subscript $0$ means that the computation are made with respect to the
metric $\dd \sigma_0^2$. 

We use the subscript $i$ when the computation are made with respect to
$g_i$ in $\overline{\Ome_i}$. Let us define the following quantity :
$$
I_n=\sum_i\int_{\Ome_i}\Div_i(\phi_n^2u\nabla_i u)\dd v_i
$$
We notice that, since $\phi_n$ is compactly supported, $I_n$ is well
defined.

In fact, because $u\Delta_i u\ge0$, we have:
\begin{align*}
I_n&=\sum_i\int_{\Ome_i}2\phi_n u\la\nabla_i \phi_n,\nabla_i u\ra_i\dd v_i
+\int_{\Ome_i}\phi_n^2\|\nabla_i u\|_i^2\dd v_i+ \int_{\Ome_i} \phi_n^2
u\Delta_i u\dd v_i\\
&\ge \sum_i\int_{\Ome_i}2\phi_n u\la\nabla_i \phi_n,\nabla_i u\ra_i\dd v_i
+\int_{\Ome_i}\phi_n^2\|\nabla_i u\|_i^2\dd v_i
\end{align*}

Because of Section~\ref{rappel}, we also have :
\begin{align*}
I_n&=\sum_i\int_{\partial\Ome_i}\phi_n^2u\la\nabla_iu,\nu_i\ra_i\dd
\boH^1_i\\
&=\frac{1}{2}\sum_{(i,j)}\left(\int_{\Gamma_{i,j}}\phi_n^2u\la\nabla_iu,
\nu_i\ra_i\dd
\boH^1_i+\int_{\Gamma_{i,j}}\phi_n^2u\la\nabla_ju,\nu_j\ra_j\dd
\boH^1_j\right)
\end{align*}
where $\nu_i$ is the outgoing normal from $\Ome_i$
along $\Gamma_{i,j}$. We notice that the results of
Subsection~\ref{selfinter}
are applied for a Riemannian metric however this Stokes formula can be
easily deduced from the Euclidean one. Let $C_{i,j}$ be the part of
$\partial O_i\cap\partial O_j$ such that $\pi(C_{i,j})=\Gamma_{i,j}$. Let
$n_i$ be the unit outgoing normal from $O_i$ in $S''$ or $\Sigma_{\eps_1}$.
We then have: 
\begin{align*}
\int_{\Gamma_{i,j}}\phi_n^2u\la\nabla_iu, \nu_i\ra_i\dd \boH^1_i
&+\int_{\Gamma_{i,j}}\phi_n^2u\la\nabla_ju,\nu_j\ra_j\dd \boH^1_j\\
&=\int_{C_{i,j}}\phi_n^2(f_K\circ \d-\eps_1)\la\nabla
(f_K\circ \d),n_i\ra\dd \boH_{\dd s^2}\\
&\quad\quad+\int_{C_{i,j}}\phi_n^2(f_K\circ
\d-\eps_1)\la\nabla (f_K\circ \d),n_j\ra\dd \boH_{\dd s^2}\\
&=\int_{C_{i,j}}\phi_n^2(f_K\circ \d-\eps_1)(f_K'\circ\d)\la\nabla
\d,n_i+n_j\ra\dd \boH_{\dd s^2}
\end{align*}
where $\phi_n$ is extended to $M_+(\eps_1)$ by $\phi_n(p)=\phi_n(\pi(p))$.

By construction, a point $p\in S_0$ is such that $\d(p)\le \d(q)$ for any
$q\in \pi^{-1}(\pi(p))\cap(S''\cup\Sigma_{\eps_1})$ this implies that
along $C_{i,j}$, $\la\nabla \d,n_i+n_j\ra\ge 0$. Hence since
$\phi_n^2(f_K\circ \d-\eps_1)(f_K'\circ\d)\le 0$, we obtain $I_n\le 0$.
This proves that
$$
\sum_i\int_{\Ome_i}2\phi_n u\la\nabla_i \phi_n,\nabla_i u\ra_i\dd v_i
+\int_{\Ome_i}\phi_n^2\|\nabla_i u\|_i^2\dd v_i\le 0.
$$
Thus:
\begin{align*}
\sum_i\int_{\Ome_i}\phi_n^2\|\nabla_i u\|_i^2\dd v_i &\le
-2\sum_i\int_{\Ome_i}\phi_n u\la\nabla_i \phi_n,\nabla_i u\ra_i\dd v_i\\
&\le 2\left(\sum_i\int_{\Ome_i}\phi_n^2\|\nabla_i u\|_i^2\dd
v_i\right)^\frac{1}{2}
\left(\sum_i\int_{\Ome_i}u^2\|\nabla_i\phi_n\|_i^2\dd
v_i\right)^\frac{1}{2}.
\end{align*}

Thus 
$$
\sum_i\int_{\Ome_i}\phi_n^2\|\nabla_i u\|_i^2\dd v_i \le 4
\sum_i\int_{\Ome_i}u^2\|\nabla_i\phi_n\|_i^2\dd v_i.
$$

The function $u$ is bounded and the metric $g_i$ and $\dd\sigma_0^2$ are
$k$-quasi-isometric so there exists a constant $C$ which does not
depend on $i$ and $n$ such that 
$$
\int_{\Ome_i}u^2\|\nabla_i\phi_n\|_i^2\dd v_i\le
C\int_{\Ome_i}\|\nabla_0\phi_n\|_0^2\dd v_0.
$$
Hence :
$$
\sum_i\int_{\Ome_i\cap\phi_n^{-1}(1)}\|\nabla_i u\|_i^2\dd v_i \le
4C\int_{\Sigma_0}\|\nabla_0\phi_n\|_0^2\dd v_0.
$$
Taking the limit $n\rightarrow +\infty$ we obtain :
$$
\sum_i\int_{\Ome_i}\|\nabla_i u\|_i^2\dd v_i =0.
$$
This implies that $u$ is constant so $S_0\subset\Sigma_{\eps_1}$, this
gives the contradiction we look for and Theorem~\ref{halfspace} is proved.

\subsubsection{$S\looparrowright M_+(\eps)$ and $H_0\ge0$}

In the second case, the only difference is the construction of the
superharmonic function. It is in fact simpler since we do not have
to control $\la n,\xi\ra$. From \eqref{surharmonic2}, we have
$$
\Delta_{S''}f_K\circ \d \le \left(f_K''(\d)+f_K'(\d)(
\kappa_1\cos^2\theta+ \kappa_2\sin^2 \theta)\right)\sin^2\phi
$$
There is still a constant $C$ such that $\max(|\kappa_1|,|\kappa_2 |) \le
C$. Then considering $K=C$ we get 
\begin{align*}
\Delta_{S''}f_K\circ \d&\le
\left(f_K''(\d)+f_K'(\d)C\right)\sin^2\phi\\
&\le 0
\end{align*}
$f_K\circ \d$ is then superharmonic and this gives also a contradiction.

\subsection{Hypotheses and stable surfaces}

In this subsection, we want to make a remark about the hypothesis of
Theorem~\ref{halfspace}. 

Let $\Sigma_0$ be as in the theorem and assume that $\Sigma_0$ has
constant mean
curvature $H_0$. Applying the Jacobi operator to the constant function
$1$, the hypothesis about the mean curvature of the equidistant surfaces
implies that $0\ge L(1)=-(2Ric(n,n)+|A|^2)$ along $\Sigma_0$.

Now assume that $\Sigma_0$ is stable, since $\Sigma_0$ is parabolic there
exists a sequence of compactly supported smooth functions $(\phi_k)_k$ such
that $0\le \phi_k\le 1$, $(\phi_k^{-1}(1))_k$ is a compact exhaustion of
$\Sigma_0$ and  
$$
\lim_k\int_{\Sigma_0}\|\nabla\phi_k\|^2=0
$$

Then by stability we get:
\begin{align*}
0\ge\int_{\phi_k^{-1}(1)}-(2Ric(n,n)+|A|^2)&\ge
\int_{\Sigma_0}-(2Ric(n,n)+|A|^2)\phi_k^2\\
&\ge\int_{\Sigma_0}\phi_kL\phi_k-\int_{\Sigma_0}\|\nabla\phi_k\|^2\\
&\ge-\int_{\Sigma_0}\|\nabla\phi_k\|^2.
\end{align*}

Taking the limit as $k$ goes to $+\infty$, we obtain $2Ric(n,n)+|A|^2=0$
along $\Sigma_0$. This implies that, at first order, the equidistant
surfaces to $\Sigma_0$ have constant mean curvature $H_0$.

Now if the equidistant surfaces have constant mean curvature $H_0$, we get
$0=L(1)$ and $2Ric(n,n)+|A|^2=0$. $\Sigma_0$ is then a stable cmc $H_0$
surface.

If $\Sigma_0$ is not stable, we see that there exists $\eps'>0$ such that
no $\Sigma_t$, $0<t<\eps'$, has constant mean curvature $H_0$. Thus
Theorem~\ref{halfspace} says that there is no constant mean curvature
$H_0$ surface in $M_\pm(\eps')$ (with good  orientation in $M_+(\eps')$).

%% file: fighalfspace3.pstex_t
\begin{picture}(0,0)%
\includegraphics{fighalfspace3.pstex}%
\end{picture}%
\setlength{\unitlength}{4144sp}%
\begingroup\makeatletter\ifx\SetFigFontNFSS\undefined%
\gdef\SetFigFontNFSS#1#2#3#4#5{%
  \reset@font\fontsize{#1}{#2pt}%
  \fontfamily{#3}\fontseries{#4}\fontshape{#5}%
  \selectfont}%
\fi\endgroup%
\begin{picture}(5697,2394)(413,-2375)
\put(4431,-188){\makebox(0,0)[lb]{\smash{{\SetFigFontNFSS{14}{16.8}{\rmdefault}{\mddefault}{\updefault}{\color[rgb]{0,0,0}$\Sigma_{\eps_1}$}%
}}}}
\put(5391,-1126){\makebox(0,0)[lb]{\smash{{\SetFigFontNFSS{14}{16.8}{\rmdefault}{\mddefault}{\updefault}{\color[rgb]{0,0,0}$S''$}%
}}}}
\put(5585,-1868){\makebox(0,0)[lb]{\smash{{\SetFigFontNFSS{14}{16.8}{\rmdefault}{\mddefault}{\updefault}{\color[rgb]{0,0,0}$S_0$}%
}}}}
\put(2008,-1463){\makebox(0,0)[lb]{\smash{{\SetFigFontNFSS{14}{16.8}{\rmdefault}{\mddefault}{\updefault}{\color[rgb]{0,0,0}$D_0$}%
}}}}
\put(3268,-2288){\makebox(0,0)[lb]{\smash{{\SetFigFontNFSS{14}{16.8}{\rmdefault}{\mddefault}{\updefault}{\color[rgb]{0,0,0}$\Sigma_0$}%
}}}}
\end{picture}%

%% file: halfspace6.tex
In this section, we prove a halfspace result when the ambient space is a
Lie group with a left invariant Riemannian metric.

Let $G$ be a $3$-dimensional connected Lie group and $F$ be a normal
properly embedded $2$-dimensional Lie subgroup. We denote by $\gotg$ and
$\gotf$ the associated Lie algebras.

Let $\dd s^2$ be a left invariant metric on $G$. $F$ is then a constant
mean curvature surface in $G$. Do we have a halfspace theorem with respect
to $F$? In fact for any $g\in G$, the coset $gF$ is also a constant mean
curvature surface in $G$. Since the left multiplication by $g$ is
an isometry, the halfspace problem is the same as the one for $F$.

Let $X\in\gotg$ the left invariant unit vector field which is normal to $F$
at $e$. Let $Y$ be a left invariant vector field, we have 
$$
\la\nabla_XX,Y\ra=-\la[X,Y],X\ra
$$
Since $F$ is normal, for any $Y\in\gotf$, $[X,Y]\in\gotf$. Then $X$ normal
to $\gotf$ implies that $\nabla_XX=0$. Then $t\mapsto \exp(tX)$ is the
geodesic from $e$ with speed $X$ at $e$.

The map $F\times\R\rightarrow G,\ (f,t)\mapsto f\exp(tX)$ is onto. Let
$t_0>0$ be the lower bound of $\{t>0|\,\exp(tX)\in F\}$. If $t_0$ exists,
$F$ does not separate $G$ and the above map is bijective on
$F\times[0,t_0)$. If $t_0=+\infty$, $G$ is diffeomorphic to $F\times\R$ and
$F$ separates $G$.

We have the following halfspace result.

\begin{prop}\label{liegroup}
Let $G$ be a $3$-dimensional connected Lie group with a left invariant
metric $\dd s^2$. Let $F$ be a normal properly embedded $2$-dimensional Lie
subgroup of $G$ which is parabolic for the left invariant metric. We
denote by $H_0$ the mean curvature of $F$. Let $S$ be a properly immersed
constant mean curvature $H_0$ surface in $G$ with no boundary.
\begin{itemize}
\item If $F$ does not separate $G$ and $S$ is included in $G\setminus
F$, $S$ is a coset $gF$.
\item If $F$ separates $G$ and $S$ is included in the mean convex side of
$F$, $S$ is a coset $gF$.
\item If $F$ separates $G$ and $S$ is included in the non mean convex side
of $F$ and is well oriented with respect to $F$, $S$ is a coset $gF$.
\end{itemize}
\end{prop}

Let us just explain what is well oriented with respect to $F$. If
$G_+$ is the non mean convex side of $F$ and $D$ is the connected
component of $G_+\setminus S$ containing $F$, we ask that along
$S\cap \partial D$ the mean curvature vector of $S$ points into $D$.

\begin{proof}
Let $X\in\gotg$ still denote the left invariant unit vector field which
is normal to $F$ at $e$ and points into the mean convex side. Let
$s\mapsto g(s)=\exp(sX)$ be the geodesic curve from the unit element $e\in
G$ normal to $\gotf=T_eF$. For any $f\in F$, $s\mapsto fg(s)$ is the
geodesic curve from $f\in F$ normal to $T_fF$. So the equidistant to $F$ at
distance $t$ is $Fg(t)$. Since $F$ is normal $Fg(t)=g(t)F$, thus the
equidistant to $t$ has the same mean curvature as $F$ and the norm of
its second fundamental form is constant. We denote by $F_t$ this
equidistant. Depending on the case, $G$ can be parametrized by
$F\times[0,t_0)$ or $F\times\R$ such that $F\times\{s\}$ is an equidistant
surface to $F$. The mean convex side is the part included in $F\times\R_+$
(there is a change of sign with respect to the preceding section). The
projection map $\pi_s$ from $F_s$ to $F_0$ is given by the right
multiplication by $g(s)^{-1}$.

Let $s_0\in\R$ be such that $F\times[0,s_0]\cap S$ is non empty.
$F\times[0,s_0]$ is then a outside or inside regular $s_0$-neighborhood
that satisfies the hypothesis about the mean curvature of the equidistant
($F\times[0,s_0]$ is regular because the right multiplication by
$g(s)^{-1}$ is quasi-isometric).
Moreover, $F$ is parabolic, so Theorem~\ref{halfspace} applies and $S$ is
an equidistant surface to $F$ \textit{i.e.} a coset $gF$.
\end{proof}

Proposition~\ref{liegroup} can be applied in several situations.
\begin{example}
For $(\mu,\nu)\in\R^2$, let us define $G(\mu,\nu)=\R^3$ with the following
Lie group structure:
$$
(x_1,y_1,z_1)\cdot(x_2,y_2,z_2)=(x_1+e^{\mu z_1}x_2,x_2+e^{\nu
z_1}y_2,z_1+z_2)
$$
and the left invariant metric:
$$
g_{\mu,\nu}=e^{-2\mu z}\dd x^2+e^{-2\nu z}\dd y^2+\dd z^2.
$$
The surfaces $\{z=t\}$ are the cosets of the normal Lie subgroup $\{z=0\}$
which is parabolic and has constant mean curvature
$H(\mu,\nu)=|\mu+\nu|/2$. Thus we have a halfspace result for cmc
$H(\mu,\nu)$ surfaces in $G(\mu,\nu)$ with respect to these horizontal
surfaces $\{z=t\}$.
\begin{prop}\label{application}
Let $S$ be a properly immersed constant mean curvature $H(\mu,\nu)$
surface in $G(\mu,\nu)$ with no boundary.
\begin{enumerate}
\item If $S$ is included in the mean convex side of one $\{z=t\}$, $S$
is equal to one $\{z=t'\}$.
\item If $S$ is included in the non mean convex side of one $\{z=t\}$ and
$S$ is well oriented with respect to it, $S$ is equal to one
$\{z=t'\}$.
\end{enumerate}
\end{prop}
Actually, some special cases of this proposition are already known. When
$\mu=0=\nu$, $G(0,0)$ is just $\R^3$ with its Euclidean flat metric and
Proposition~\ref{application} is the classical halfspace theorem for
minimal surfaces with respect to a plane \cite{HoMe}. When $\mu=\nu=c\neq
0$, $G(c,c)$ is the hyperbolic space $\H^3$ and we recover the halfspace
theorem with respect to horospheres \cite{RodRos}. When $(\mu,\nu)=(0,c)$
($c\neq 0$), $G(0,c)$ is $\H^2(c)\times\R$ and $\Sigma_t$ are the
horocylinders so we get the halfspace theorem proved in \cite{HaRoSp}. The
last special case is when $\mu+\nu=0$, $G(-c,c)$ is the model space
$\mathrm{Sol}_3$ and we obtain the halfspace theorem given by Theorem 1.5
in \cite{DaMeRo}.
\end{example}
\begin{example}
Let $G$ be the Heisenberg group $\Nil_3$ \textit{i.e.} $G=\R^3$ with the
Lie group structure:
$$
(x_1,y_1,z_1)\cdot(x_2,y_2,z_2)=(x_1+x_2,y_1+y_2,z_1+z_2+\frac{x_1y_2}{2}
-\frac{x_2y_1}{2})
$$
and the left invariant metric:
$$
\dd s^2=\dd x^2+\dd y^2+(\frac{1}{2}(y\dd x-x\dd y)+\dd z)^2
$$
For any $(a,b)\in\R^2\setminus\{(0,0)\}$, the surfaces $\{ax+by=c\}$ are
the cosets of the normal Lie subgroup $\{ax+by=0\}$ which is parabolic and
has mean curvature $0$. Thus we recover the halfspace theorem for
``vertical minimal planes'' proved by B.~Daniel and L.~Hauswirth in
\cite{DaHa}.
\end{example}

\begin{example}
Let us consider, on $G=\R\times\R_+^*\times\R$, the Lie group structure:
$$
(x_1,y_1,z_1)\cdot(x_2,y_2,z_2)=(x_1+y_1x_2,y_1y_2,z_1+z_2)
$$
with the left invariant metric:
$$
\dd s^2=\frac{1}{y^2}(\dd x^2+\dd y^2)+(\frac{1}{y}\dd x+\dd z)^2.
$$
For any $t$, the surfaces $\Sigma_t=\{y=e^t\}$ are the cosets of the
normal Lie subgroup $\{y=1\}$ which is parabolic and has mean curvature
vector $\dis\frac{1}{2}\der{}{y}$. Thus Proposition~\ref{liegroup} can be
apply to obtain Proposition~\ref{psl2}.

In fact, as Riemannian manifolds, $G$ is isometric to
$\widetilde{PSL}_2(\R)$ that is $\R\times\R_+^*\times\R$ with the same
Riemannain metric but with a different Lie group structure (the expression
of the Lie group structure is not easy to write so we prefer to omit it).
For this Lie group structure, the metric is still left invariant but
$\{y=1\}$ is no longer a subgroup. However we obtain the following
halfspace result.
\begin{prop}\label{psl2}
Let $S$ be a properly immersed constant mean curvature
$\frac{1}{2}$ surface in $\widetilde{PSL}_2(\R)$ with no boundary.
\begin{enumerate}
\item If $S$ is included in the mean convex side of one $\Sigma_t$, $S$
is equal to one $\Sigma_{t'}$.
\item If $S$ is included in the non mean convex side of one $\Sigma_t$
and $S$ is well oriented with respect to it, $S$ is equal to one
$\Sigma_{t'}$.
\end{enumerate}
\end{prop}
In fact, the projection map $(x,y,z)\mapsto (x,y)$ is a
Riemannian submersion from $\widetilde{PSL}_2(\R)$ to $\H^2$. So
the surfaces $\Sigma_t$ that foliate $\widetilde{PSL}_2(\R)$ are called
``vertical horocylinders'' since they are the fiber over horocycles in
$\H^2$. Proposition~\ref{psl2} is then a halfspace result with respect to
the vertical horocylinders in $\widetilde{PSL}_2(\R)$.

The author recently learns that this result is also proved by
Carlos Penafiel in \cite{Pen}
\end{example}